\documentclass[11pt,reqno]{amsart} 
\usepackage[utf8]{inputenc}
\usepackage{amsmath}
\usepackage{amscd}
\usepackage{amssymb, bm}
\usepackage{amsthm}
\usepackage{mathrsfs}
\usepackage{xspace}
\usepackage[all,tips]{xy}
\usepackage[dvips]{graphicx}
\usepackage{syntonly}
\usepackage{hyperref}
\usepackage{graphics, fullpage,color, epsfig,url}
\usepackage{indentfirst}
\usepackage{esint}
\usepackage{enumitem}
\usepackage[dvipsnames]{xcolor}
\usepackage{tensor}
\usepackage{amsrefs, eucal}
\usepackage{bbm}
\usepackage{ulem}
\numberwithin{equation}{section}

\newcommand{\be}{\begin{equation}}
\newcommand{\ee}{\end{equation}}

\newcommand{\vs}{\vspace{0.2cm}}

\newcommand{\ben}{\begin{equation*}}
\newcommand{\een}{\end{equation*}}

\theoremstyle{plain}
\newtheorem{theorem}{Theorem}[section]
\newtheorem{Remark}[theorem]{Remark}
\newtheorem{Definition}[theorem]{Definition}
\newtheorem{Proposition}[theorem]{Proposition}

\newtheorem{Lemma}[theorem]{Lemma}

\DeclareMathOperator{\R}{\mathbb{R}}

\numberwithin{equation}{section}

\title{A monotonicity formula for a  semilinear fractional  parabolic equation}
\author{Ignacio Bustamante}
\address{IMERL, Facultad de Ingeniería, Universidad de la República,
Montevideo, Uruguay}
\email{{\href{mailto:}{ibustamante@fing.edu.uy}}}

\begin{document}

\begin{abstract}
   We apply a high-dimensional parabolic-to-elliptic transformation to the fully fractional, semilinear heat equation $(\partial_t -\Delta)^s u = |u|^{p-1}u$, where $0<s<1$ and $p>1$, and establish a monotonicity formula for the corresponding extension problem. This result serves as a fractional analogue of the Giga-Kohn monotonicity formula for the local equation $\partial_t u - \Delta u = |u|^{p-1}u.$ The method of proof may be applicable to other nonlinear and nonlocal settings.
\end{abstract}

\maketitle 
\section{Introduction}

Many central problems in geometric analysis and partial differential equations rely on understanding how certain quantities vary under scaling, deformation, or evolution with respect to time. A powerful technique, applied for example in Almgren’s study of harmonic functions \cite{A83}, Alt-Caffarelli-Friedman’s analysis of free boundaries \cite{Alt1984} or Huisken's works on the mean curvature flow  \cite{H90}, involves constructing integral functionals that exhibit monotonic behavior along a natural parameter. These monotonicity formulas have since become indispensable tools, offering compactness, regularity, and rigidity results in settings ranging from elliptic PDEs to geometric flows \cite{CM12}. 

\vs

These methods have also been successfully extended to nonlocal operators. For example, for the fractional Laplacian introduced by Riesz \cite{Riesz},
$$(-\Delta)^s v(x) := \frac{4^s\Gamma(d/2+s)}{\pi^{d/2}|\Gamma(-s)|}\lim_{r \to 0^+}\int_{\R^d \setminus B_r(x)} \frac{v(x)-v(z)}{|x-z|^{d+2s}}dz, \quad \mbox{ where } 0<s<1,$$ 
several monotonicity formulas have been established. Among these are an Almgren frequency-type monotonicity formula due to Caffarelli and Silvestre \cite{CS} and an Alt-Caffarelli-Friedman type monotonicity formula proved by Terracini, Verzini and Zilio \cite{Terracini2016}.
These results typically rely on the celebrated Caffarelli-Silvestre extension \cite{CS}, which interprets the fractional Laplacian as a Dirichlet-to-Neumann operator for a degenerate but local PDE on the half space $\R^{d+1}_+$.

\vs

For the parabolic counterpart of the fractional Laplacian, the fractional heat operator defined in \cite{Riesz} as
\be\label{definition fractional heat}
(\partial_t-\Delta)^su(x,t):=\int_{-\infty}^t \int_{\R^d}(u(x,t) - u(z,\tau)) K_s(x-z,t-\tau)dzd\tau, \quad \mbox{ where } 0<s<1,
\ee
 and
\be\label{kernel K}
K_s(z,\tau) := \frac{1}{(4\pi)^{d/2}|\Gamma(-s)|}\frac{e^{-|z|^2/4\tau} }{\tau^{d/2+1+s}},\ee 
monotonic quantities have also been identified. In particular, an Almgren frequency-type parabolic monotonicity formula was established by Stinga and Torrea \cite{Stinga2017}, and an Alt-Caffarelli-Friedman type parabolic monotonicity formula was recently proved by Davey and Smit Vega Garcia \cite{DaveySmit25}. As in the local case, clear parallels exist between the elliptic and parabolic formulas, with the latter typically depending on carefully constructed backward solutions to heat-type equations.  
Parabolic monotonicty formulas also rely on an extension problem, developed by   Nystr\"{o}m and  Sande \cite{Nystrm2016} and independently by Stinga and Torrea \cite{Stinga2017}.
Deriving  parabolic monotonic quantities, whether in local or nonlocal settings, is often challenging, particularly in nonlinear contexts.

 \vs

This article focuses on the semilinear equation
\be\label{Fractional lane parabolic}
(\partial_t -\Delta)^s u = |u|^{p-1}u,
\ee
where $0<s<1$ and $p>1$. The local case ($s=1$) is well-understood, with established results for its well-posedness, regularity theory, and blow-up profiles (see \cite{Quittner2019} and references therein). For the local case, Giga and Kohn \cite{Giga1985} derived a fundamental monotonicity formula, which plays a crucial role in characterizing the blow-up profiles of solutions.

\vs

Several problems related to (\ref{Fractional lane parabolic}) have recently appeared in the literature: Chen, Guo and Li considered equations of the form $(\partial_t -\Delta)^su =f(x,t)$ where $f\geq 0$, and established regularity results \cite{ChenGuoli}. Ma, Guo and Zhang  analyzed the case where  $f=f(u)$ on $B_1(0)\times \R$ \cite{Ma2024}. For problems with spatial dependence, Chen and Guo  proved nonexistence of bounded solutions when $f(x,u) = x_1u^p$ \cite{Guochen}. Banerjee and Garofalo  examined the linear case $f(x,t,u) = V(x,t)u$, obtaining unique continuation results \cite{BG18, BG24}. The work most closely related to ours is that of Ferreira and de Pablo, who studied the case where $f(x,t,u)=u^p$ on $(0,T)$, with nonnegative memory data, i.e., imposing $u(x,t) =g(x,t)$ for $t \leq 0$ \cite{Ferreira2024}. Note that a maximum principle is available for this operator (see Theorem 1.5 of \cite{Stinga2017}) ensuring that nonnegative memory data yields nonnegative solutions to \eqref{Fractional lane parabolic}, and therefore, equation (\ref{Fractional lane parabolic}) coincides with the problem studied in \cite{Ferreira2024} when considering appropriate memory data.

\vs

The stationary solutions to (\ref{Fractional lane parabolic}) correspond to those of the fractional Lane-Emden equation
\be\label{lane emden}
(-\Delta)^s u = |u|^{p-1}u,
\ee
since, as shown in \cite{Stinga2017},  although the fractional heat operator is nonlocal in both space and time, it reduces to the fractional Laplacian $(-\Delta)^s$ when applied to a function that solely depends on $x$,
\ben
(\partial_t - \Delta)^s u(x) = (-\Delta)^s u(x).
\een
Equation (\ref{lane emden}) has been extensively studied, with many classical results extended to the nonlocal setting (see \cite{Li04, Chen05, Dvila2017} and references therein). 
Of particular relevance to this work is the monotonicity formula found by Dávila, Dupaigne and Wei \cite{Dvila2017}, which was used to classify solutions of finite Morse index. This quantity, which we discuss in Section \ref{section 2}, can be viewed as the fractional analogue of the local monotonicity formula for the equation $-\Delta u=|u|^{p-1}u$, as discussed by Fazly and Shahgholian \cite{FazlyShahgholian} (see also the article by Pacard \cite{Pacard1993} for a similar monotonicity formula in the case $-\Delta u = u^{p}$). 

\vs

In this article, we derive a new monotonicity formula for solutions of (\ref{Fractional lane parabolic}), which can be interpreted as the parabolic counterpart of the Dávila-Dupaigne-Wei monotonicity formula for the fractional Lane-Emden equation. Specifically, let $u=u(x,t)$ be a backward solution to (\ref{Fractional lane parabolic}) on a time interval $(0,T_I)$, where $T_I>0$, that is, $u(x,t):=\bar{u}(x,-t)$ where $\bar{u}$  is a solution of \eqref{Fractional lane parabolic} (see Section \ref{ssec:fractional_heat}). Due to the nonlocal nature of (\ref{definition fractional heat}), $u$ must be defined in $(0,+\infty)$, so we may either prescribe $u(\cdot,t) = f(\cdot,t)$ for $t \geq T_I$, or consider backward ancient solutions of  (\ref{Fractional lane parabolic}) instead. Let $U(x_0,x,t)$ be its corresponding extension, defined in $\R_+ \times \R^d \times (0,T_I)$. Then, under appropriate growth and regularity assumptions, the function
\begin{equation}\label{J}
    \begin{split}
    \mathcal{J}(t) :=& \int_{\R^{d+1}_+} x_0^{1-2s}\frac{|\nabla U|^2}{2} \mathcal{G}_s\, dX - \frac{\eta_{s}}{p+1} \int_{\R^d} |u|^{p+1}\tilde{\mathcal{G}}_s\,dx + \frac{s}{p-1}\int_{\R^{d+1}_+} x_0^{1-2s} \frac{U^2}{2t} \mathcal{G}_s\,dX,
\end{split}
\end{equation}
where $X=(x_0,x) \in \R^{d+1}_+ :=\R_+ \times \R^{d}$, is non-decreasing for the time-reversed variable $t$, and its derivative is explicitly given by,
\be\label{DJ}
\begin{split}
\frac{d}{dt}\mathcal{J}(t) &= \int_{\R^{d+1}_+} x_0^{1-2s}\left( \partial_t U +\frac{X}{2t} .\nabla U + \frac{2 s}{p-1}\frac{U}{2t} \right)^2\mathcal{G}_s\,dX,
\end{split}
\ee
see Theorem \ref{main}. Here,
\ben
     \mathcal{G}_s(X,t) := t^{\frac{2s}{p-1}+1}\mathcal{G}(X,t) \quad  \mbox{ and }\quad 
    \tilde{\mathcal{G}}_s (x,t) := t^{\frac{2s}{p-1}+1}\mathcal{G}((0,x),t),
    \een
    are appropriate rescalings of the fundamental solution $\mathcal{G}$ for the extension problem of the equation $(\partial_t -\Delta)^s u = 0,$
\be\label{kernel}
\mathcal{G}(X,t) = \frac{1}{(4\pi)^{d/2}\Gamma(s)}\frac{e^{-|X|^2/4t} }{t^{d/2+1-s}},\ee
where $X \in \R^{d+1}_+$,  $t>0,$  and $\eta_s$ is a constant given by,
    \be\label{eta s}
    \eta_s := \frac{2s|\Gamma(-s)|}{4^s\Gamma(s)}.
    \ee 
    
    This result provides a fractional analogue of the Giga-Kohn monotonicity formula \cite{Giga1985} (compare with the presentation in \cite{Ecker2005}). To our knowledge, it constitutes the first monotonicity formula established for nonlinear, fully fractional heat equations. Moreover, it offers a new tool to analyze singularities for the problem studied by Ferreira-de Pablo \cite{Ferreira2024}.  
    
    The growth and regularity assumptions hold, for example, if $u$ is twice differentiable and its first and second derivatives are bounded, although these restrictions arise from the method of proof itself, and could potentially be relaxed through more direct approaches.

\vs
The proof we present relies on a high-dimensional parabolic-to-elliptic transformation for the extension problem of backward solutions of \eqref{Fractional lane parabolic}. By applying this transformation, we recast the parabolic extension problem into a high-dimensional elliptic problem, which is shown to hold if and only if the parabolic one does. The elliptic problem we obtain structurally resembles the extension problem for the fractional Lane-Emden equation \eqref{lane emden}, though additional terms arise. Nevertheless, the resemblance allows for an adaptation of the existing monotonicity results for equation (\ref{lane emden}). By adapting the monotonicity formula, we obtain an almost monotonic quantity, in the sense that it becomes monotonic in the high-dimensional limit. Finally, through careful analysis of this limit, we derive the parabolic monotonicity formula for the original problem, and recover its derivative.

\vs
The strategy of proof we employ here was first explored by Perelman  in the local setting, by means of which the celebrated reduced volume for the Ricci flow was derived \cite{Per02}, and has attracted increasing attention in recent years. We now provide a brief account of these developments. Tao \cite{T09l} and Svérak \cite{Sve11} discussed this approach, with Svérak emphasizing its probabilistic interpretation. Building on the probabilistic viewpoint, Davey \cite{Dav18} derived several well-known parabolic monotonicity formulas from their elliptic counterparts, thereby recovering a number of classical results. Davey and Smit Vega Garcia \cite{DS24} further extended this method to variable coefficients and established a new Alt–Caffarelli–Friedman monotonicity formula in that setting. Recent work by Reiris and the author \cite{BR} has also shown that Perelman's entropy also arises as a high-dimensional limit of an elliptic monotonicity formula, Colding's monotonic volume, when applied to Perelman's original high-dimensional space, showing that both reduced volume and entropy can be understood as originating from a unified, high-dimensional Ricci-flat framework. 
More recently, and during the final stages of preparation of this work, an article by Davey and Smit Vega Garcia appeared on arXiv, addressing the adaptation of these techniques to fractional operators in the linear case $(\partial_t - \Delta)^s u = 0$ \cite{DaveySmit25}. 

It is important to remark that, in contrast to the probabilistic approach adopted by Davey  in \cite{Dav18} and Davey and Smit Vega Garcia in \cite{DS24, DaveySmit25}, our method follows the geometric perspective originally proposed by Perelman, which proves effective in treating the nonlinear case by controlling the volume element restricted to the boundary of the higher-dimensional space. Moreover, this geometric approach also enables us to obtain not only a monotonic quantity, but also the explicit formula for its derivative.

\vs

Finally, we would like to note that the same techniques can indeed derive the classical Giga-Kohn monotonicity formula for $\partial_t u -\Delta u = |u|^{p-1}u$  from the monotonicity formula  for  $-\Delta u = |u|^{p-1}u$ discussed by Fazly and Shahgholian \cite{FazlyShahgholian}.  Since the computations are carried out in a similar manner, we omit them in the interest of brevity. The methods we present here hold potential for application in other nonlinear settings, as well as for systems of equations. 

\vs

\textbf{Organization of the article.} The article is structured as follows: Section \ref{section 2} presents the relevant properties  of the fractional Laplacian and fractional heat operator, including their Caffarelli-Silvestre extension problems. We also briefly discuss the fractional Lane-Emden equation and its monotonicity formula. Then, we develop the parabolic-to-elliptic transformation and formulate the corresponding non-homogeneous elliptic extension problem. Section \ref{section 3} establishes an almost-monotonicity formula for solutions of the non-homogeneous extension problem. In Section \ref{section 4}, we derive explicit volume element formulas for relevant subsets and hypersurfaces in the high-dimensional space and prove the convergence results needed for our main result, Theorem \ref{main}. The proof of the theorem is completed in Section \ref{section 5}.

\section{Fractional operators and parabolic-to-elliptic transformations}\label{section 2}
In this section, we review some results concerning the fractional Laplacian, the fractional Lane-Emden equation, and the fractional heat operator. We then derive the non-homogeneous elliptic problem solved by the backward extension of the parabolic equation we address. For convenience, we will work on $N$-dimensional spaces for fractional elliptic problems, and on $d$-dimensional spaces for fractional parabolic ones.
\subsection{The fractional Laplacian and fractional Lane-Emden equation.}
Let $0<s<1$, and let $v:\R^N \to \R$  such that  $v \in C^{2s+\varepsilon}(\R^N)$ for some $\varepsilon>0$. Assume
 $$\int_{\R^N} \frac{|v(z)|}{(1+|z|)^{N+2s}}\,dz < +\infty,$$
 so that its fractional Laplacian, $(-\Delta)^sv$, is well defined. From now on, we denote $a:=1-2s$.

 \vs

A celebrated result by Caffarelli and Silvestre \cite{CS}, shows that the fractional Laplacian may be interpreted as a Dirichlet-to-Neumann operator for a local, degenerate operator on the half-space $\R^{N+1}_+$. Specifically, let $(z_0,z)\in \R^{N+1}_+$, where we assume $z_0 \in \R_+$. If $v \in C^{2s+\varepsilon}(\R^N) \cap L^1(\R^N, (1+|z|)^{-N-2s}dz),$ define
$$
V(z_0,z) := \int_{\R^{N}} v(z-y)P(z_0,y)\, dy,
$$
where $P$ is the Poisson kernel,
$$P(z_0,z) := C_{N,a}z_0^{1-a}|(z_0,z)|^{-(N+1-a)}, $$
and $C_{N,a}>0$ is chosen so that $\int_{\R^N} P(z_0,z)dz=1$. Then $V \in C^{2}(\R^{N+1}_+) \cap C(\overline{\R^{N+1}_+}),$ $z^{a}_0\partial_{z_0} V \in C(\overline{\R^{N+1}_+})$ and $V$ is a solution of the extension problem,
\be\label{elliptic extension}
\begin{cases}
\nabla \cdot (z_0^{a}\nabla V) = 0\quad &\mbox{ for } (z_0,z) \in \R^{N+1}_+, \\
V(0,z) = v(z) \quad &\mbox{ for } z \in \R^{N}.\\
\end{cases}
\ee
Moreover, $V$ obeys,
\ben
-\lim_{z_0 \to 0^+}z_0^a \partial_{z_0}V(z_0,z)  = \kappa_s (-\Delta)^s v,
\een
with
$$\kappa_s := \frac{\Gamma(1-s)}{2^{2s-1}\Gamma(s)}.$$
Suppose now that $v$ is a solution of the fractional Lane-Emden equation (\ref{lane emden}). Then,
\be\label{derivative lane emden eq}
-\lim_{z_0 \to 0^+}z_0^a \partial_{z_0}V(z_0,z) = \kappa_s |v|^{p-1}v\, (z),
\ee
as discussed in \cite{Dvila2017}. 
Observe that the first equation in (\ref{elliptic extension}) is equivalent to
\be\label{first eq extension lane}
\Delta_z V + \frac{a}{z_0}\partial_{z_0}V +\partial^2_{z_0}V = 0 \, \mbox{ for } (z_0,z) \in \R^{N+1}_+.
\ee
For such $V$, the following monotonicity formula is known.
\begin{theorem}[Theorem 1.4 in \cite{Dvila2017}] 
 Let  $V(z_0,z) \in C^{2}(\mathbb{R}_{+}^{N+1}) \cap C(\overline{\mathbb{R}_{+}^{N+1}})$, such that $V$ obeys  $($\ref{elliptic extension}$)$ and $($\ref{derivative lane emden eq}$)$, and suppose $ z_0^a \partial_{z_0} V \in C(\overline{\mathbb{R}_{+}^{N+1}}) $. For $ R>0$, let
\be\label{homogeneous monotonicity}
\begin{aligned}
 E\left(R\right) :=&R^{2 s \frac{p+1}{p-1}-N}\left(\frac{1}{2} \int_{\mathbb{R}_{+}^{N+1} \cap B^{N+1}_{R}} z_0^a|\nabla V|^{2} d z d z_0-\frac{\kappa_{s}}{p+1} \int_{\partial \mathbb{R}_{+}^{N+1} \cap B^{N+1}_R}|V|^{p+1} \,d z\right) \\
& \quad +R^{2 s \frac{p+1}{p-1}-N-1} \frac{s}{p-1} \int_{\partial B^{N+1}_R \cap \mathbb{R}_{+}^{N+1}} z_0^a V^{2}\, d \sigma.
\end{aligned}
\ee
Then, $E$ is a non-decreasing function of $R$. Moreover,
\be
\frac{d E}{d R}=R^{2 s \frac{p+1}{p-1}-N} \int_{\partial B^{N+1}_R \cap \mathbb{R}_{+}^{N+1}} z_0^{a}\left(\frac{\partial V}{\partial r}+\frac{2 s}{p-1} \frac{V}{r}\right)^{2} \,d \sigma.
\ee
\end{theorem}

Here, $B^{N+1}_R$ denotes the Euclidean ball in $\mathbb{R}^{N+1}$ centered at the origin of radius $R,$ $ \sigma$ is the $N$-dimensional Hausdorff measure restricted to the hypersurface $\partial B^{N+1}_R$, $r=|(z_0,z)|$ the Euclidean norm of a point $(z_0,z) \in \mathbb{R}_{+}^{N+1}$, and $\partial_{r}= \frac{(z_0,z)}{r} \cdot \nabla$ is the radial derivative.

\subsection{The fractional heat operator.}\label{ssec:fractional_heat}

Now, we turn our attention to some relevant results concerning the fractional heat operator. Let $0<s<1$, and assume that $u$ solves
\be\label{Parabolic LE on interval}
(\partial_t -\Delta)^s u = |u|^{p-1}u \quad \mbox{for } (x,t) \in \R^d \times (-T_I,T_F),
\ee
where $T_I,T_F >0.$
By definition, for the fractional heat operator (\ref{definition fractional heat}) to exist at a time $t_0 \in (-T_I,T_F)$, it must also exist for every time $t \leq t_0 $. Moreover, 
$(\partial_t -\Delta)^su$ is well defined for $t \in (-T_I,T_F)$ provided $u$ belongs to the parabolic Hölder continuous space $u \in C_{t,x}^{2s+\varepsilon}(\R^d \times (-\infty, T_F))$ for some $\varepsilon>0$ \cite{Stinga2017}. We now recall this definition.
 
 \begin{Definition}
     We say  $u: \Omega \to \R$ is in the parabolic Hölder continuous space $u \in C_{t,x}^{\gamma}(\Omega)$ with $0<\gamma\leq 1$ if $u \in L^{\infty}(\Omega)$, and there exists $\gamma>0$ such that
     $$|u(x,t) -u(z,s)| \leq C (|x-z|^2 +|t-s|)^{\gamma/2},$$
      for every $(x,t), (z,s) \in \Omega$. 
      
For $1< \gamma \leq 2$, we say that $u \in C^{\gamma}_{t,x}(\Omega)$ if $u \in L^{\infty}(\Omega)$, $u$ is $\gamma/2$-Hölder continuous in $t$ uniformly in $x$ and its gradient $\nabla_xu$ is $(\gamma-1)-$Hölder continuous in $x$ uniformly in $t$.
 \end{Definition}

Let $u:(-\infty,T_F) \to \R$ be a parabolic Hölder continuous function of order $2s+\varepsilon$ such that $u$ solves (\ref{Parabolic LE on interval}). The data $u|_{(-\infty,-T_I]}$ may either be prescribed (this is referred to as the memory problem in the literature and some well-posedness results are available for nonnegative memory data, see \cite{Ferreira2024}), or we may instead consider ancient solutions, that is, solutions of
$
(\partial_t -\Delta)^s u = |u|^{p-1}u$ for $(x,t) \in \R^d \times (-\infty,T_F).$
\vs

Now, define the parabolic extension $U$ of the function $u$ as
\be\label{parabolic extension}
U(x_0,x,t) := \int_{0}^{\infty} \int_{\R^d} P^s_{x_0} (z,\tau) u(x-z,t-\tau) \, dz d\tau,
\ee
where $X=(x_0,x) \in \R_+ \times \R^d$ and $P_{x_0}^s(z,\tau)$ is the fractional Poisson kernel,
$$P_{x_0}^s(z,\tau) = \frac{1}{4^{d/2+s} \pi^{d/2}\Gamma(s) } \frac{x_0^{2s}}{\tau^{d/2+1+s}}e^{-(x_0^2 +|z|^2)/4\tau}, $$
where the constant is chosen so that \be\label{Property Kernel}
\int_0^{\infty} \int_{\R^d} P_{x_0}^s(z,\tau)dzd\tau =1.
\ee
As shown in \cite{Nystrm2016,Stinga2017},  $U$ is well defined whenever $u \in C_{t,x}^{2s+\varepsilon}(\R^d \times (-\infty, T_F))$  and  satisfies two key properties: first, it solves the extension problem,
\be\label{fractional parabolic}
\begin{cases}\partial_{t} U=\Delta_x U+\frac{a}{x_0} \partial_{x_0}U+\partial_{x_0}^2U, & \text { for }(X,t) \in \mathbb{R}^{d+1}_+\times (-T_I,T_F), \\ U(0,x,t )=u(x,t),&\text { for }(x,t) \in \mathbb{R}^{d}\times (-T_I,T_F),
\end{cases}
\ee
and second, we can recover the fractional heat operator using the normal derivative at the boundary $\partial\R^{d+1}_+,$
$$
\eta_s |u|^{p-1}u=\eta_s (\partial_t -\Delta)^s u=- \lim _{x_0 \to 0^{+}} x_0^{a} \partial_{x_0} U(x_0,x,t),
$$
where $\eta_s$ is the constant defined in (\ref{eta s}).
The proof of this fact presented in \cite{Stinga2017} relies on the Fourier transform. Nevertheless, as discussed in  \cite[Section 2]{Ferreira2024}, the extension $U$ obeys these properties whenever the integrals involved are well defined.

\vs

The extension (\ref{parabolic extension}) can alternatively be expressed in terms of the fundamental solution $\mathcal{G}$. If $u$ is a solution of the master equation $$(\partial_t - \Delta)^s u = h,$$ for some regular enough $h$, the solution of the parabolic extension problem $($\ref{fractional parabolic}$)$ can be written as,
$$\label{explicit U}
U(x_0,x,t) := \int_0^{\infty} \int_{\mathbb{R}^d} \mathcal{G}(x_0,z,\tau)h(x-z,t-\tau)\,dzd\tau,
$$
where $\mathcal{G}$ is defined in $(\ref{kernel})$, see \cite{Stinga2017}. As before, $U$ obeys
$$
\eta_s h(x,t)=- \lim _{x_0 \to 0^{+}} x_0^{1-2 s} \partial_{x_0} U(x_0,x,t).
$$
A similar expression can also be proved for solutions of the memory problem of equation (\ref{Fractional lane parabolic}), if we assume the memory data is regular and decays fast enough as $t \to -\infty$, see \cite{Ferreira2024}.
Moreover, it can be checked that the function $\mathcal{G}$ obeys
$$\lim_{x_0 \to 0^+} x_0^a \partial_{x_0}\mathcal{G}(x_0,x,t) = 0,$$
for any $(x,t) \in \R^d \times \R_+$, and $$\partial_t \mathcal{G}= \Delta_x \mathcal{G} + \frac{a}{x_0}\partial_{x_0} \mathcal{G} + \partial^2_{x_0} \mathcal{G},$$ for any positive time $t$.

\vs

 Now, consider backward solutions of equation (\ref{Fractional lane parabolic}), which are defined as follows: for any function $g:\R^d \times \R$, we let $\bar{g}(x,t) := g(x,-t)$  denote its time reversal. Then, given a function $u:(-\infty, T_F) \to \R$  which solves (\ref{Parabolic LE on interval}), we have that $\bar{u}: \R^d \times (-T_F,+\infty)$ solves  
\be\label{back frac}
(-\partial_t -\Delta)^s u = |u|^{p-1}u \, \mbox{ for } (x,t) \in \R^d \times (-T_F,T_I),
\ee
 where
\ben
\begin{split}
(-\partial_t -\Delta)^s u(x,t)&=
(\partial_t-\Delta)^s\bar{u}(x,-t) \\
&=\int_t^{\infty} \int_{\R^d}(u(x,t) - u(z,\tau)) \overline{K_s}(x-z,t-\tau)\,dzd\tau,
\end{split}
\een
see \cite{Stinga2017, DaveySmit25} for a detailed discussion.In order to see this, notice that if $u$ solves (\ref{Parabolic LE on interval}), then
\be\label{backwards lane}
(-\partial_t -\Delta)^s \bar{u}(x,t) =  (\partial_t -\Delta)^s u(x,-t) = |u|^{p-1}u(x,-t)=|\bar{u}|^{p-1}\bar{u}(x,t).
\ee
Similarly, starting from a function  $u \in C_{t,x}^{2s+\varepsilon}(\R^d \times (-T_F,+\infty))$ which solves the backward fractional heat equation (\ref{back frac}), its time reversal $\bar{u}$ is a solution of the forward fractional heat equation (\ref{Parabolic LE on interval}), and from (\ref{parabolic extension}) it follows that the extension associated to $u$ is 
\be\label{backwards extension parabolic}
U(x_0,x,t) = \bar{U}(x_0,x,-t) = \int_{-\infty}^{0} \int_{\R^d} \overline{P^s_{x_0}} (z,\tau) u(x-z,t-\tau)\, dz d\tau,
\ee
where $\bar{U}$ is the extension (\ref{parabolic extension}) associated to  $\bar{u}$. By (\ref{fractional parabolic}), $U$ is a solution to the backward extension problem, 
\be\label{backwards fractional parabolic}
\begin{cases}\partial_{t} U+\Delta_x U+\frac{a}{x_0} \partial_{x_0}U+\partial_{x_0}^2U=0, & \text { for }(X,t) \in \mathbb{R}^{d+1}_+\times (-T_F,T_I), \\  U(0,x,t )=u(x,t),&\text { for }(x,t) \in \mathbb{R}^{d}\times (-T_F,T_I). \end{cases}
\ee
As before, we find that 
\be\label{derivative fractional}
-\lim _{x_0 \to 0^{+}} x_0^{a} \partial_{x_0} U(x_0,x,t) =\eta_s (-\partial_t -\Delta)^s u(x,t)= \eta_s |u|^{p-1}u(x,t),
\ee
for every $(x,t) \in \R^d \times (-T_F,T_I).$

\vs

From this point forward, we will exclusively consider backward solutions,  defined in $\R^d \times (-T_F, +\infty)$  for some $T_F>0$, and obey equation (\ref{backwards lane}) in $(-T_F, T_I)$. We will restrict our attention to compact time intervals which, by a time translation, we assume to be $[0,T]$ for some $T\in (0,T_I)$. Backward solutions will be denoted by $u$, and their extensions will be denoted by $U$. We will work within the class of functions we now define.

\begin{Definition}\label{class C}
    We say that $U: \R^{d+1}_+ \times [0,T] \to \R$ belongs to the function class $\mathcal{U}([0,T])$
    if $U \in C^2(\R^{d+1}_+\times [0,T])\cap C(\overline{\R^{d+1}_+} \times [0,T])$ and the following hold:
    \begin{enumerate}[label=$($\alph*$)$]
    \item  $\lim\limits _{x_0 \to 0^+}x_0^{1+a}\partial_{t}U(x_0,x,t)=0$ for $(x,t) \in \R^{d} \times [0,T]$ and $x_0^a \partial_{x_0}U \in C(\overline{\R^{d+1}_+} \times [0,T]).$
        \item The functions $f_U$ belong to $ C((0,T);L^2(\R^{d+1}_+,\, x_0^a dX))$, and
        $$\sup_{t \in (0,T)} ||f_U(X,t)||^2_{L^2(\R^{d+1}_+,\, x_0^a dX)} < +\infty,$$
        for any of the following $f_U$:
     $$ e^{-|X|^2/8t}\nabla U, \quad   e^{-|X|^2/8t}U , \quad e^{-|X|^2/8t} X\cdot\nabla U, \quad e^{-|X|^2/8t}\partial_t U.$$
    \item We have 
    $$ \sup_{t \in (0,T)} \big|\big|
            g_U(X,t) e^{-|X|^2/4t} \big|\big|_{L^1(\R^{d+1}_+ , \,x_0^adX)}< +\infty,$$
        for any of the following $g_U$:
        $$  UH, \quad \partial_tU H,  \quad H( X\cdot\nabla U), $$
        where \be\label{H}H:=2(X,t)\cdot \nabla_{(X,t)}\partial_t U + (d+1+a)\partial_t U.
        \ee
        Here $\nabla$ denotes the gradient with respect to the $X$ variables, and $\nabla_{(X,t)}$ denotes the gradient with respect to the $(X,t)$ variables.
\end{enumerate}
\end{Definition}
    Items $(b)$ and $(c)$ of Definition \ref{class C} impose moderate growth controls over $U$ and its derivatives. While these conditions suffice for our proof, the results may remain valid under less restrictive assumptions. We also note that related classes of functions (adapted to the method of proof employed and the hypotheses needed for each problem) have been considered when deriving monotonicity formulae for variable coefficient parabolic operators in \cite{DS24} and, recently, for solutions of the extension problem of the fractional parabolic equation $(-\partial_t -\Delta)^su =0$  in \cite{DaveySmit25}.
The following proposition presents a class of functions for which its backward extensions obey Definition \ref{class C}.
  
   \begin{Proposition}
Assume $u:\R^d \times [0,+\infty) \to \R$ is a parabolic Hölder continuous function of order $2s+\varepsilon$ such that its time reversal $\bar{u}:\R^{d} \times (-\infty,0] \to \R$ is an ancient solution of $(\ref{Fractional lane parabolic})$. 
If $u \in C^2(\R^d \times [0,+\infty))$ and  its first and second derivatives are bounded, then $U \in \mathcal{U}([0,T])$ for any $T>0$.

Similarly, let $u: \R^d \times (-\infty, T_F) \to \R$ be a solution of the problem with memory for equation \eqref{Parabolic LE on interval}, such that the memory data is twice differentiable and satisfies $|(\partial_t-\Delta)^sf(x,t)|\leq C$. Then, its backward extension satisfies $U \in \mathcal{U}([0,T])$ for any $0<T<T_I$.

\begin{proof}
Since $|u|<C$, using the bounds on the derivatives together with expressions  $(\ref{backwards extension parabolic})$ and $(\ref{Property Kernel})$, we can show that $|U|<C$, $|\partial_{x_i}U|<C$ for $i\in \{1,\dots, d\}$, and $|\partial_t U|<C$. The fact that $|\partial_{x_0}U|< Cx_0^{-a}$ follows by using the last line of the representation formula $(1.5)$ in \cite{Stinga2017}. 

 For the problem with memory, we may use the representation formula $(2.23)$ in \cite{Ferreira2024} to prove the last bound instead. The other bounds for $U$ and its derivatives are also valid in this case, and they can be obtained using expression $(\ref{backwards extension parabolic})$ directly.   

        From the previous estimates, we see that for any fixed $t>0$, $f_U(\cdot,t) \in L^2(\R^{d+1}_+,x_0^a dX)$, the supremum of their squared norms is finite and the bounds are explicitly computable.  Continuity of the functions $f_U(\cdot,t)$ to $L^2(\R^{d+1}_+,\, x_0^a dX)$ follows from the Dominated Convergence Theorem, by employing the Gaussian decay. Condition $(c)$ follows in a similar manner. The regularity properties of $U$ discussed in \cite{Stinga2017} ensure the ones in our definition hold. 
\end{proof}
   \end{Proposition}
    
\subsection{Parabolic-to-elliptic transformations}
Given a function $U: \R^{d+1}_+\times \R \to \R$, we define its $n$-dimensional lift $V_n:\R^{d+n+1}_+ \to \R$ of $U$ as follows. Let $(z_0,z,y) \in \R_+ \times \R^d \times \R^{n}$. Then, set
\be\label{variables}
\begin{cases}
    x_0 &= \sqrt{n}z_0,\\
    x &= \sqrt{n}z,\\
    2t &= R^2 = z_0^2 + |z|^2 + |y|^2.
\end{cases}
\ee
Then, let
$$\mathcal{F}_n(z_0,z,y) := (\sqrt{n}z_0,\sqrt{n}z, R^2/2),$$
and use it to define,
\be\label{lift}
V_n (z_0,z,y) := U\circ \mathcal{F}_n(z_0,z,y) = U(x_0,x,t).
\ee
Observe that if $U$ is defined on  a region $\R^{d+1}_+ \times [0,T)$, then $V_n$ is defined on the region $ \R^{d+n+1}_+\cap B^{d+n+1}_{{\sqrt{2T}}}$. 
A direct application of the chain rule yields the following.

\begin{Lemma}\label{chain rule}
Let   $V_n: \R^{n+d+1}_+ \to \R$ as in \eqref{lift}. Then $V_n$ satisfies:
\begin{align*}
 \partial_{z_i}V_n &= \sqrt{n}\partial_{x_i}U +z_i\partial_t U, \\
 \partial_{y_j}V_n &= y_j\partial_t U, \\
 \partial_{ z_i}^2V_n &= n \partial^2_{ x_i}U +2x_i\partial^2_{x_i t} U + z_i^2\partial_{t}^2 U +\partial_tU,\\
 \partial_{ y_j }^2 V_n &= \partial_t U +y_j^2 \partial^2_{t} U, \\
 |\nabla V_n|^2 &= n|\nabla U|^2 + 2(X\cdot\nabla U)\partial_t U +2t (\partial_tU)^2,
\end{align*}
for any $i=0,1,\dots,d$, and $j = 1,\dots, n.$
In particular, 
\be
\frac{a}{z_0}\partial_{z_0} V_n = \frac{a}{z_0}\left(\sqrt{n}\partial_{x_0}U +z_0\partial_t U\right)=\frac{na}{x_0}\partial_{x_0}U + a\partial_t U,
\ee
and
\be\label{lap V_n}
\Delta_{(z,y)} V_n + \frac{a}{z_0} \partial_{z_0} V_n + \partial^2_{ z_0} V_n = n\left( \partial_t U+\Delta_x U +\frac{a}{x_0}\partial_{x_0}U+\partial^2_{x_0 }U \right)+ H,
\ee
where $H$ is defined in $($\ref{H}$)$.
\end{Lemma}
\vs

Let  $u \in C^{2s+\varepsilon}_{t,x}(\R^d \times (-T_F,+\infty))$ be a solution of the backward equation \eqref{back frac}, and let $[0,T] \subset (-T_F,T_I)$.  Let $U:\R^{d+1}_+ \times [0,T] \to \R$ be its associated extension, and assume that $U \in \mathcal{U}([0,T])$. We now establish the extension problems satisfied by the lifts $V_n$ of $U$.
First, notice that since  $U$ is a solution of (\ref{backwards fractional parabolic}), by (\ref{lap V_n}) we have,
\be
\Delta_{(z,y)} V_n + \frac{a}{z_0} \partial_{z_0} V_n + \partial^2_{z_0} V_n =H,
\ee
and $H$ does not depend on $n$.
\vs
Then, notice that
$$
\lim_{z_0 \to 0^+} V_n(z_0,z,y)=\lim_{z_0 \to 0^+}U(\sqrt{n}z_0,\sqrt{n}z, R^2/2) = \lim_{z_0 \to 0^+}U(\sqrt{n}z_0,\sqrt{n}z, (z_0^2 + |z|^2 +|y|^2)/2).
$$
Since $U \in \mathcal{U}([0,T])$, $U$ is continuous in $\overline{\R^{d+1}_+} \times [0,T]$ and therefore,
$$\lim_{z_0 \to 0^+} V_n(z_0,z,y)= U(0,\sqrt{n}z,(|z|^2+|y|^2)/2).$$

Then, using (\ref{backwards fractional parabolic}),
\be
\lim_{z_0 \to 0^+} V_n(z_0,z,y) = U(0,x,(|z|^2+|y|^2)/2) = u(x,(|z|^2+|y|^2)/2),
\ee
and, by Lemma \ref{chain rule},
\begin{equation*}
\begin{split}
z_0^a\partial_{z_0} V_n(z_0,z,y) &= z_0^a \left(\sqrt{n}\partial_{x_0} U + z_0 \partial_t U\right)\\ 
& = \left( \frac{x_0}{\sqrt{n}}\right)^a\left(\sqrt{n}\partial_{x_0}U + \frac{x_0}{\sqrt{n}}\partial_t U\right)\\
& = n^{\frac{1-a}{2}}x_0^a\left(\partial_{x_0}U + \frac{x_0}{n}\partial_t U\right).
\end{split}
\end{equation*}
Using (a) of Definition \ref{class C} and   (\ref{derivative fractional}),
\begin{equation*}
\begin{split}
\lim_{z_0 \to 0^+} z_0^a\partial_{z_0} V_n(z_0,z,y) & =  \lim_{x_0 \to 0^+} n^{\frac{1-a}{2}}x_0^a\partial_{x_0}U\left(x_0,x,t\right)\\
& = -\eta_s n^{\frac{1-a}{2}}|u|^{p-1} u(x,(|z|^2+|y|^2)/2).
\end{split}
\end{equation*}
Finally, observe that 
$$-\eta_s n^{\frac{1-a}{2}}|u|^{p-1} u(x,(|z|^2+|y|^2)/2) = -\eta_s n^{\frac{1-a}{2}}|V_n|^{p-1}V_n(0,z,y).$$
Combining the previous computations, we obtain the following.
\begin{Proposition}\label{system for V_n}
Let  $u \in C^{2s+\varepsilon}_{t,x}(\R^d \times (-T_F,+\infty))$ such that $u$ solves \eqref{back frac}, and let $[0,T] \subset (-T_F,T_I)$. 
 Let $U$ be its associated extension. If $ U \in \mathcal{U}( [0,T])$, then $V_n$ obeys
\be\label{non-homogeneous}
\begin{cases}
\nabla \cdot (z_0^{a}\nabla V_n) = z_0^a H \circ \mathcal{F}_n = n^{-a/2}x_0^aH \quad &\mbox{ in } \R^{d+n+1}_+\cap B^{d+n+1}_{{\sqrt{2T}}} \\
V_n(0,z,y) =  u(x,(|z|^2+|y|^2)/2) \quad &\mbox{ for } (z,y) \in  B^{d+n}_{{\sqrt{2T}}},\\
\end{cases}
\ee
and
\be\label{normal der}
-\lim_{z_0 \to 0^+}z_0^a \partial_{z_0}V_n = \eta_s n^{\frac{1-a}{2}}|V_n|^{p-1}V_n(0,z,y) \quad \mbox{ for } (z,y) \in \partial \R^{d+n+1}_+ \cap B^{d+n+1}_{{\sqrt{2T}}}.
\ee
Moreover, since $V_n(z_0,z,y) = U(x_0,x,t),$ $$V_n \in C^2(\R^{d+n+1}_+\cap B^{d+n+1}_{{\sqrt{2T}}}) \cap C(\overline{\R^{d+n+1}_+}\cap B^{d+n+1}_{{\sqrt{2T}}}),$$ and $$z_0^a \partial_{z_0}V_n \in C(\overline{\R^{d+n+1}_+}\cap B^{d+n+1}_{{\sqrt{2T}}}).$$
\end{Proposition}

\section{An almost monotonicity formula for $V_n$.}\label{section 3}
We now establish an ``almost monotonicity formula'' for solutions of the non-homogeneous extension problem (\ref{non-homogeneous}). This formula is considered almost monotonic in the sense that it becomes monotonic as $n \to \infty$. This fact will be evident later, since we will show that the contribution of the source terms decays to zero by relating the integral quantities defined for the high-dimensional lift with those defined for $U$ in $\R^{d+1}_+ \times [0,T]$. To prove our result, we adapt the arguments used in the proof of Theorem 1.4 in \cite{Dvila2017}.

\begin{theorem}\label{almost monotonicity}
 Let  $u \in C^{2s+\varepsilon}_{t,x}(\R^d \times (-T_F,+\infty))$ such that $u$ solves \eqref{back frac}, and let $[0,T] \subset (-T_F,T_I)$.  Let $U \in \mathcal{U}([0,T])$ be its associated extension, and $V_n$ its $n$-dimensional lift. Then the function $\mathcal{E}_n:(0,\sqrt{2T}) \to \R$ defined as,
\be\label{c_ne_n}
\begin{aligned}
 \mathcal{E}_n(R):=&R^{2 s \frac{p+1}{p-1}-N}\left(\frac{1}{2} \int_{\mathbb{R}_{+}^{N+1} \cap B^{N+1}_R} z_0^a|\nabla V_n|^{2} \,dz_0dzdy-n^{\frac{1-a}{2}}\frac{\eta_{s}}{p+1} \int_{\partial \mathbb{R}_{+}^{N+1} \cap B^{N+1}_R}|V_n|^{p+1} \, d zdy\right) \\
& \quad+R^{2 s \frac{p+1}{p-1}-N-1} \frac{s}{p-1} \int_{\partial B^{N+1}_R  \cap \mathbb{R}_{+}^{N+1}} z_0^a V_n^{2} \,d \sigma,
\end{aligned}
\ee
obeys,
\begin{equation}\label{der c_n}
\begin{split}
\frac{d \mathcal{E}_n}{d R}(R)= &R^{2 s \frac{p+1}{p-1}-N} \int_{\partial B^{N+1}_R \cap \mathbb{R}_{+}^{N+1}} z_0^{a}\left(\frac{\partial V_n}{\partial r}+\frac{2 s}{p-1} \frac{V_n}{r}\right)^{2} \,d \sigma \\
& - R^{2s\frac{p+1}{p-1}-N-1}\int_{\mathbb{R}_{+}^{N+1} \cap B_{R}^{N+1}}\left(\frac{2s}{p-1}V_n +r \frac{\partial V_n}{\partial r} \right)(z_0^{a}H \circ \mathcal{F}_n)\,dz_0dzdy,
\end{split}
\end{equation}
where $N=n+d$, $r=|(z_0,z,y)|,$ and $\partial V_n  /\partial_{r}= \frac{(z_0,z,y)}{r} \cdot \nabla V_n $ denotes the radial derivative of $V_n$.
\begin{proof}
Since $U \in \mathcal{U}([0,T]),$ Proposition \ref{system for V_n} holds. We reduce the integral to a fixed domain by introducing a rescaled function. For $(z_0,z,y) \in \mathbb{R}_{+}^{N+1}$, let
$$
W((z_0,z,y); R):=R^{\frac{2 s}{p-1}} V_n(Rz_0,Rz,Ry).
$$
A direct computation using the chain rule shows that $W$ satisfies
\ben
\begin{split}
\nabla \cdot (z_0^{a}\nabla W) ((z_0,z,y);R) & = R^{2s/(p-1)+2-a}\left(\nabla \cdot (z_0^{a}\nabla V_n)\right)(R z_0,R z, R y)\\
&=R^{2s/(p-1)+2-a}(z_0^{a}H \circ \mathcal{F}_n)(R z_0,R z,R y),
\end{split}
\een
and, for the Neumann condition,
\ben
\begin{split}
-\lim _{z_0 \to 0^+} z_0^{a} \partial_{z_0} W((z_0,z,y);R) 
&=- R^{2s/(p-1)+1}\lim _{z_0 \to 0^+} z_0^{a} \partial_{z_0} V_n(R z_0, R z, R y)\\ 
&=- R^{2s/(p-1)+1-a}\lim _{z_0 \to 0^+} (R z_0)^{a} \partial_{z_0} V_n(R z_0, R z, R y) \\
& = R^{2sp/(p-1)}\eta_sn^{\frac{1-a}{2}}(|V_n|^{p-1} V_n )(0, R z, R y) \\
& = n^{\frac{1-a}{2}}\eta_s (|W|^{p-1}W)((0, z, y); R).
\end{split}
\een
Therefore, $W$ obeys
\be\label{equation for W}
\begin{cases} 
\nabla \cdot (z_0^{a}\nabla W) ((z_0,z,y);R)=R^{2s/(p-1)+2-a}(z_0^{a}H \circ \mathcal{F}_n)(R z_0,R z,Ry),\\
-\lim\limits_{z_0 \to 0^+} (z_0^{a} \partial_{z_0} W)((z_0,z,y);R)= n^{\frac{1-a}{2}}\eta_s |W|^{p-1}W((0,z,y);R).
\end{cases}
\ee

We denote by $$\frac{\partial W}{\partial r}((z_0,z,y);R) := \frac{(z_0,z,y)}{r}\cdot \nabla W((z_0,z,y);R)$$ 
the radial derivative of $W$, and compute its derivative with respect to the parameter $R$:
\ben
\frac{\partial W}{\partial R}((z_0,z,y);R) 
= \frac{2s}{p-1}R^{2s/(p-1)-1}V_n(Rz_0,Rz,Ry) + R^{2s/(p-1)}\nabla V_n(Rz_0,Rz,Ry)\cdot(z_0,z,y).
\een
By the chain rule, $\nabla W((z_0,z,y);R) = R^{2s/(p-1)+1}\nabla V_n(Rz_0,Rz,Ry)$, and therefore
$$
r\frac{\partial W}{\partial r}((z_0,z,y);R) 
= (z_0,z,y)\cdot \nabla W((z_0,z,y);R) 
= R^{2s/(p-1)+1}(z_0,z,y)\cdot\nabla V_n(Rz_0,Rz,Ry).
$$
Combining these two expressions, we obtain the identity
\be\label{radial identity}
R\frac{\partial W}{\partial R} = \frac{2s}{p-1}W + r\frac{\partial W}{\partial r}.
\ee

Now, define
\be
\begin{split}
\tilde{\mathcal{E}}_n (V_n;R):= R^{2 s \frac{p+1}{p-1}-N}\left( \int_{\mathbb{R}_{+}^{N+1} \cap B^{N+1}_R} z_0^a\frac{|\nabla V_n|^{2}}{2}\, d y dz d z_0
-\frac{n^{\frac{1-a}{2}}\eta_{s}}{p+1} \int_{\partial \mathbb{R}_{+}^{N+1} \cap B^{N+1}_R}|V_n|^{p+1} \,d ydz\right).
\end{split}
\ee
A scaling argument shows 
$$\tilde{\mathcal{E}}_n(V_n;R) = \tilde{\mathcal{E}}_n(W;1),$$
and differentiating the right-hand side with respect to $R$,
\be\label{dR tilde E_n}
\frac{d \tilde{\mathcal{E}}_n}{d R}(V_n ; R)
=\int_{\mathbb{R}_{+}^{N+1} \cap B_{1}^{N+1}} z_0^{a} \nabla W \cdot \nabla \frac{\partial W}{\partial R}\,d y dz d z_0
- n^{\frac{1-a}{2}}\eta_{s} \int_{\partial \mathbb{R}_{+}^{N+1} \cap B_{1}^{N+1}}|W|^{p-1}W \frac{\partial W}{\partial R}\, d y dz.
\ee
To handle the first integral on the right-hand side, we use \eqref{equation for W} to write
\be\label{boundary term E_n}
\begin{split}
z_0^{a}\nabla W \cdot \nabla \frac{\partial W}{\partial R}
&= \nabla \cdot \left(z_0^a \frac{\partial W}{\partial R} \nabla W\right) - \frac{\partial W}{\partial R} \nabla \cdot (z_0^a\nabla W)\\
&= \nabla \cdot \left(z_0^a \frac{\partial W}{\partial R} \nabla W\right) 
- \frac{\partial W}{\partial R}((z_0,z,y);R) \cdot R^{2s/(p-1)+2-a}(z_0^{a}H \circ \mathcal{F}_n)(Rz_0,Rz,Ry),
\end{split}
\ee
and to simplify notation, define
$$A((z_0,z,y);R):=R^{2s/(p-1)+2-a}(z_0^{a}H \circ \mathcal{F}_n)(Rz_0,Rz,Ry).$$
Integrating the divergence term in \eqref{boundary term E_n} over $\mathbb{R}_{+}^{N+1} \cap B_{1}^{N+1}$
yields a boundary integral over $\partial B_1^{N+1} \cap \mathbb{R}_+^{N+1}$ and a boundary integral
over $\partial\mathbb{R}_+^{N+1} \cap B_1^{N+1}$. The latter, combined with the Neumann condition
in \eqref{equation for W}, cancels exactly with the second term in \eqref{dR tilde E_n}, giving
\be\label{derivative non-homogeneous}
\frac{d \tilde{\mathcal{E}}_n}{d R}(V_n ; R) 
= \int_{\partial B_{1}^{N+1} \cap \R^{N+1}_+} z_0^{a} \frac{\partial W}{\partial r} \frac{\partial W}{\partial R}\, d \sigma 
- \int_{\mathbb{R}_{+}^{N+1} \cap B_{1}^{N+1}} \frac{\partial W}{\partial R} \, A\, dz_0d y dz.
\ee
We now use the identity \eqref{radial identity} 
on $\partial B_1^{N+1}$, where $r=1$, to obtain
\ben
\begin{aligned}
\int_{\partial B_{1}^{N+1} \cap \R^{N+1}_+} z_0^{a} \frac{\partial W}{\partial r} \frac{\partial W}{\partial R}\, d \sigma
&=R \int_{\partial B_{1}^{N+1} \cap \mathbb{R}_{+}^{N+1}} z_0^{a} \left(\frac{\partial W}{\partial R}\right)^{2}\, d \sigma
-\frac{2 s}{p-1} \int_{\partial B_{1}^{N+1} \cap \mathbb{R}_{+}^{N+1}} z_0^{a} W \frac{\partial W}{\partial R}\, d \sigma,\\
&=R \int_{\partial B_{1}^{N+1} \cap \mathbb{R}_{+}^{N+1}} z_0^{a} \left(\frac{\partial W}{\partial R}\right)^{2}\, d \sigma
-\frac{s}{p-1}\frac{d}{dR}\left(\int_{\partial B_{1}^{N+1} \cap \mathbb{R}_{+}^{N+1}} z_0^{a} W^{2} \,d \sigma\right).
\end{aligned}
\een
Substituting in \eqref{derivative non-homogeneous},
\be\label{derivative E_n in terms of V}
\begin{split}
\frac{d \tilde{\mathcal{E}}_n}{d R}(V_n ; R) 
&=R \int_{\partial B_{1}^{N+1} \cap \mathbb{R}_{+}^{N+1}} z_0^{a} \left(\frac{\partial W}{\partial R}\right)^{2}\, d \sigma
-\frac{s}{p-1}\frac{d}{dR}\left(\int_{\partial B_{1}^{N+1} \cap \mathbb{R}_{+}^{N+1}} z_0^{a} W^{2} \,d \sigma\right) \\
& \phantom{=}\,
- \int_{\mathbb{R}_{+}^{N+1} \cap B_{1}^{N+1}} \frac{\partial W}{\partial R} \, A \, dz_0dzdy.
\end{split}
\ee

It remains to scale back to $V_n$ on $B_R^{N+1}$. The first term on the right-hand side of \eqref{derivative E_n in terms of V} gives
$$
R \int_{\partial B_{1}^{N+1} \cap \mathbb{R}_{+}^{N+1}} z_0^{a} \left(\frac{\partial W}{\partial R}\right)^{2}\, d \sigma
= R^{2s\frac{p+1}{p-1}-N}\int_{\partial B_{R}^{N+1} \cap \mathbb{R}_{+}^{N+1}} z_0^{a} \left(\frac{\partial V_n}{\partial r}+\frac{2s}{p-1}\frac{V_n}{r}\right)^{2}\, d \sigma,
$$
where we used \eqref{radial identity}. For the total $R$-derivative term we obtain 
$$
\frac{s}{p-1}\frac{d}{dR}\left(\int_{\partial B_{1}^{N+1} \cap \mathbb{R}_{+}^{N+1}} z_0^{a} W^{2} \,d \sigma\right)
= \frac{s}{p-1}\frac{d}{dR}\left( R^{2s\frac{p+1}{p-1}-N-1}\int_{\partial B_{R}^{N+1} \cap \mathbb{R}_{+}^{N+1}} z_0^{a} V_n^{2} \,d \sigma\right).
$$
For the last term, using the expression for $\partial W/\partial R$ and the change of variables $(z_0,z,y) \mapsto R^{-1}(z_0,z,y)$,
\ben
\begin{split}
&\int_{\mathbb{R}_{+}^{N+1} \cap B_{1}^{N+1}} \frac{\partial W}{\partial R}((z_0,z,y);R) \, A((z_0,z,y);R)\, dz_0dzdy \\
&\phantom{==} = R^{4s/(p-1)+1-a}
\int_{\mathbb{R}_{+}^{N+1} \cap B_{1}^{N+1}}
\left[ \Bigl( \frac{2s}{p-1}V_n + (z_0,z,y)\cdot\nabla V_n\Bigr)\times  (z_0^aH\circ\mathcal F_n)
\right](Rz_0,Rz,Ry)\,dz_0dzdy\\
&\phantom{==} = R^{2s\frac{p+1}{p-1}-N-1}\int_{\mathbb{R}_{+}^{N+1} \cap B_{R}^{N+1}}
\left(\frac{2s}{p-1}V_n + r\frac{\partial V_n}{\partial r} \right)
(z_0^{a}H \circ \mathcal{F}_n)\,dz_0dzdy,
\end{split}
\een
where we used $(z_0,z,y)\cdot\nabla V_n = r\frac{\partial V_n}{\partial r}$. From here, \eqref{der c_n} follows by substituting the three scaled-back expressions above into \eqref{derivative E_n in terms of V}, and using that
\ben
\mathcal{E}_n(R) = \tilde{\mathcal{E}}_n(V_n;R) + R^{2 s \frac{p+1}{p-1}-N-1} \frac{s}{p-1} \int_{\partial B^{N+1}_R  \cap \mathbb{R}_{+}^{N+1}} z_0^a V_n^{2} \,d \sigma.
\een
\end{proof}
\end{theorem}

\section{Volume elements and convergence lemmas}\label{section 4}

Lemma \ref{chain rule} shows that the integrands appearing in Proposition \ref{almost monotonicity} can be represented in terms of the functions $U$ and $u$. To effectively express these integrals in terms of the radial variable on the high-dimensional space and
 $X=(x_0,x) \in \R^{d+1}_+$, we first examine the volume form on the subsets where integration takes place.

\subsection{Induced volume elements.}
We now let $(z_0,z,y) \in \R^{d+1} \times \R^n$. Notice that here $z_0 \in \R$. We will restrict these variables to different subsets of interest later. Now let $(l,\phi) \in \R_+ \times \mathbb{S}_1^d$ denote the polar coordinates in the $\R^{d+1}$ factor, i.e., $(z_0,z)=(l,\phi). $ Similarly, for the $\R^n$ factor, we denote the polar coordinates $\R_+ \times \mathbb{S}^{n-1}_1$ by $y = (s,\theta).$ 

\vs

The Euclidean metric in $\R^{n+d+1}$,  $g_E= dz_0^2 +dz^2+ dy^2$, can be written in these coordinates as,
\be\label{g_e}
g_E = dl^2 + l^2 d\Omega_{d}^2 + ds^2 + s^2 d\Omega_{n-1}^2,
\ee
where $d\Omega_{m}^2$ denotes the standard metric in $\mathbb{S}^{m}_1$. Notice that if we define $r=|(z_0,z,y)|$, then $s=(r^2-l^2)^{1/2}$.
We introduce coordinates $(r,X,\theta)$ in $\R^{n+d+1}$ via the map 
\be
F(r,X,\theta) = \left( \frac{X}{\sqrt{n}}, \sqrt{\left(r^2-\frac{|X|^2}{n}\right)} \, ,\theta \right),
\ee
where $r \in \R_+$, $X \in \R^{d+1},$ $\theta \in \mathbb{S}^{n-1}_1.$ Here, $(z_0,z)=X/\sqrt{n}$ and $l^2 = |X|^2/n$, ensuring $r^2 \geq l^2$ and thus $s \geq 0$. Then, the following relations hold:
\be
\begin{cases}
(l, \phi) &=(z_0,z) = X/\sqrt{n},\\ 
(s, \theta) & = y = (\sqrt{\left(r^2-|X|^2/n\right)}, \theta).
\end{cases}
\ee
Since  $r^2 = s^2 + l^2,$ on the $(n+d)$-dimensional spheres $\{r=R=\mbox{const}\}\subset \R^{n+d+1}$, we have $l^2 \leq R^2$ and
$$ds\big|_{T\{r=R\} }= -\frac{l}{s}dl\big|_{T\{r=R\} },$$
where ${T\{r=R\} }$ denotes the tangent space to the $(n+d)$-dimensional sphere $\{r=R\}.$ The metric induced by (\ref{g_e}) is,
\be
g_E\vert_{\{r=R\}} = \left(1+\frac{l^2}{R^2-l^2}\right)dl^2 + l^2 d\Omega_{d}^2 + (R^2-l^2) d\Omega_{n-1}^2.
\ee
Therefore, 
\be
\begin{split}
\sqrt{g_E}\vert_{\{r=R\}} &= \sqrt{\left(1+\frac{l^2}{r^2-l^2}\right)l^{2d}(r^2-l^2)^{n-1}}\\
& =rl^{d}(r^2-l^2)^{(n-2)/2}\\
&=r^{n-1}l^{d}\left(1- \frac{l^2}{r^2} \right)^{(n-2)/2},
\end{split}
\ee
where $l^2 \leq r^2=R^2.$ After computing the induced volume density on $\{r=R\}$, we simply replace $R$ by the ambient radial variable $r$ and append the factor $dr$ to obtain the full Euclidean volume form. Thus, the Euclidean volume form is,
$$dV = r^{n-1}l^{d}\left(1- \frac{l^2}{r^2} \right)^{(n-2)/2}dl\wedge d\phi \wedge dr \wedge d\theta.$$
Since the Euclidean metric in $\R^{d+1}$ is given by $\tilde{g}_E = dz_0^2+dz^2 = dl^2 + l^2 d\Omega_{d}^2,$ we have,
$$dz_0\wedge dz_1 \wedge \dots \wedge dz_d= l^{d}dl \wedge d\phi.$$
Now, using that 
$X=(x_0,x) = \sqrt{n}(z_0,z),$ it is straightforward to show,
$$dz_0\wedge dz_1 \wedge \dots \wedge dz_d = n^{-\frac{d+1}{2}}dx_0\wedge dx,$$
where $dx= dx_1 \wedge \dots \wedge dx_d.$
We can combine this with $l^2 =|X|^2/n \leq r^2$ to obtain,
\be
\begin{split}
dz_0 \wedge &dz \wedge dy =  n^{-(d+1)/2} r^{n-1}\left(1- \frac{|X|^2}{nr^2} \right)^{(n-2)/2}\, dx_0 \wedge dx \wedge dr \wedge d\theta.
\end{split}
\ee
Recall that $\partial \R^{n+d+1} = \{z_0 =0\} =\{x_0 =0\}$, and denote by $B^{m}_{R}$ the ball of radius $R$ in $\R^{m}$.  We now write $\partial \R^{n+d+1}_+ \cap B^{n+d+1}_{R}$ in the coordinates $(r,X,\theta)$ as,
\be\label{preimage plane}
\begin{split}
\partial \R^{n+d+1}_+ \cap B^{n+d+1}_{R} =\{(r,X,\theta) : 0 \leq r \leq R, \,  X = (x_0,x), \, x_0=0, \, x \in B^d_{\sqrt{nr^2}}, \, \theta \in \mathbb{S}^{n-1}_1 \},
\end{split}
\ee
and since $\partial_{z_0}=\sqrt{n}\partial_{x_0}$,  the induced volume form over this region is given by,
\ben
\begin{split}dz\wedge dy&= \iota_{\partial_{z_0}}dV\lvert_{\{z_0=0\}\cap B^{n+d+1}_{R}}\\ 
& = n^{-d/2}r^{n-1}\left(1- \frac{|x|^2}{nr^2} \right)^{(n-2)/2}\, dx \wedge  dr \wedge d\theta.
\end{split}
\een
 Similarly, we obtain
\be\label{preimage boundary}
\partial B_{R}^{n+d+1} \cap \R^{n+d+1}_+ =\{(r,X,\theta) : r=R,\,  X \in \R^{d+1}_+ \cap B^{d+1}_{\sqrt{nR^2}},\,  \theta \in \mathbb{S}^{n-1}_1 \}.
\ee
The induced volume form on this region is,
\ben
\begin{split}
d\sigma&=\iota_{\partial_r}dV\lvert_{\{r=R\}} \\ &=   n^{-(d+1)/2}R^{n-1}\left(1- \frac{|X|^2}{nR^2} \right)^{(n-2)/2}\,dx_0\wedge dx \wedge d\theta.
\end{split}
\een
Finally, 
\be\label{preimage ball}
B_{R}^{n+d+1} \cap \R^{n+d+1}_+ = 
\{(r,X,\theta) : 0 \leq r \leq R,\,  X \in \R^{d+1}_+ \cap B^{d+1}_{\sqrt{nr^2}}, \, \theta \in \mathbb{S}^{n-1}_1 \},
\ee
and the volume form on this region is,
\ben
\begin{split}
dz_0 \wedge  d z \wedge dy =  n^{-(d+1)/2} r^{n-1}\left(1- \frac{|X|^2}{nr^2} \right)^{(n-2)/2} dx_0 \wedge dx \wedge dr \wedge d\theta.
\end{split}
\een
Define $G_n : \R^{d+1}\times \R_+ \to \R$ by
$$G_n(X,t) := \left(1- \frac{|X|^2}{2nt} \right)^{(n-2)/2}\chi_{\R^{d+1}_+ \cap B^{d+1}_{\sqrt{2nt}}}(X),$$
and $\tilde{G}_n : \R^{d}\times \R_+ \to \R$ by
$$\tilde{G}_n(x,t) : = \left(1- \frac{|x|^2}{2nt} \right)^{(n-2)/2}\chi_{B^{d}_{\sqrt{2nt}}}(x).$$
Using $G_n$ and $\tilde{G}_n$, by (\ref{preimage plane}), (\ref{preimage boundary}) and (\ref{preimage ball}) we can write the resulting volume elements as,
\be\label{volume in plane}
dzdy\vert_{\partial \R^{n+d+1}_+ \cap B_R} =  n^{-d/2}r^{n-1}\tilde{G}_n(x,r^2/2)dxdrd\theta,
\ee
\be\label{volume in H_R}
d\sigma\vert_{\partial B_{R}^{n+d+1} \cap \R^{n+d+1}_+}  =n^{-(d+1)/2}R^{n-1}G_n(X,R^2/2)dXd\theta,
\ee
and
\be\label{volume in ball}
dz_0dzdy\vert_{B_{R}^{n+d+1} \cap \R^{n+d+1}_+} = \ n^{-(d+1)/2} r^{n-1}G_n(X,r^2/2)dXdrd\theta,
\ee
where $dX=dx_0dx.$
\subsection{Convergence lemmas.}

Now, we establish some convergence lemmas for the integral quantities defined from $V_n$. Our analysis begins with an  observation about the limits of the functions $G_n$ and $\tilde{G}_n$.

Define
\be
G(X,t) := e^{-|X|^2/4t} \quad \mbox{ and } \quad \tilde{G}(x,t) := e^{-|x|^2/4t}, 
\ee
where $G$ is a function on $\R^{d+1}_+\times \R_+$,and $\tilde{G}$ is a function on $\R^{d}\times \R_+$.
As shown in Lemma 3.1 of \cite{DS24}, we have
$$
G_n(X,t) \to G(X,t) \quad \mbox{ and } \quad \tilde{G}_n(x,t) \to \tilde{G}(x,t),
$$
uniformly in $\R^{d+1}_+ \times \{t\geq t_0\}$ and $\R^{d} \times \{t\geq t_0\}$ respectively, for any $t_0>0$. Moreover, it is also proved there that
 $$ G_n(X,t) < C_d G(X,t),$$
for $(X,t) \in \R^{d+1}_+\times \R_+$ and, similarly we have,
  $$ \tilde{G}_n(x,t) < C_d G(x,t),$$
for some constant $C_d>0$ (in fact, for our definition of $G_n$ and $G$, it is straightforward to show that  $G_n \leq eG$, and $\tilde{G}_n \leq e\tilde{G}$). The Dominated Convergence Theorem now ensures that, if $f(\cdot,t)$ is integrable with respect to $G(\cdot,t)\, x_0^adX$, it is also integrable with respect to $G_n(\cdot,t)\,x_0^adX$ for every $n \in \mathbb{N}, t>0$, and similarly for the measures $\tilde{G}(\cdot, t)\, dx$ and $\tilde{G}_n(\cdot, t) \,dx$ respectively. We have the following lemma.

\begin{Lemma}\label{Lemma convergence hn}
    Let $0<t_0<t_1$. Suppose  $f:\R^{d+1}_+ \times [t_0,t_1]\to \R$ is continuous and such that $fG \in C([t_0,t_1]; L^1(\R^{d+1}_+,\, x_0^adX))$. Define
    $$h_n(t) : = \int_{\R^{d+1}_+} fG_n(X,t)\,x_0^adX,$$
    and 
    $$h(t) : = \int_{\R^{d+1}_+} fG(X,t)\,x_0^adX.$$
    Then $h_n \to h$ uniformly in $[t_0,t_1]$.
\end{Lemma}
The proof of Lemma \ref{Lemma convergence hn} can be found in the Appendix.
\begin{Remark}\label{Lemma for up}
By examining the proof of Lemma \ref{Lemma convergence hn}, we see that the result also holds if $\tilde{f}:\R^{d} \times [t_0,t_1]\to \R$ is continuous and $\tilde{f}\tilde{G} \in C([t_0,t_1]; L^1(\R^{d}, dx))$, where the functions $h_n$ and $h$ are replaced by
    $\tilde{h}_n(t) : = \int_{\R^{d}} \tilde{f}(x,t)\tilde{G}_n(x,t)\,dx$  and $\tilde{h}(t) : = \int_{\R^{d}} \tilde{f}(x,t)\tilde{G}(x,t)\,dx$ respectively. 

In particular, if $u \in C^{2s+\varepsilon}_{t,x}(\R^d \times (-T_F,+\infty))$, this holds for $\tilde{f}=|u|^{p+1}$ on compact intervals of $(0,T)$. To show this, first observe that $|u|^{p+1}\tilde{G} <C\tilde{G}$ for some positive constant $C$. Continuity of $|u|^{p+1}\tilde{G}(\cdot,t)$ in $L^1(\R^d,dx)$ follows  from the continuity of $u$ and the Dominated Convergence Theorem applied to sequences $t_n \to t$. Moreover, $|u|^{p+1}\tilde{G} <C\tilde{G}$ also implies that
$$\sup_{t \in (0,T)} \int_{\R^d} |u|^{p+1}\tilde{G}(x,t) \,dx < \infty.$$
We will use both properties in the following section.    
\end{Remark}

\begin{Remark}\label{lemma for derivative term}
    If $U \in \mathcal{U}([0,T])$, then Lemma \ref{Lemma convergence hn} applies to the functions 
    $f(X,t) = |\nabla U|^2$, $U^2$, $|X \cdot \nabla U|^2$, and $(\partial_t U)^2$ 
    on compact subintervals of $(0,T)$.
    Indeed, by property (b) of Definition \ref{class C}, the functions
    \[
        e^{-|X|^2/8t}\nabla U,\quad 
        e^{-|X|^2/8t}U,\quad 
        e^{-|X|^2/8t} X\cdot\nabla U,\quad 
        e^{-|X|^2/8t}\partial_t U,
    \]
    belong to $C((0,T);L^2(\R^{d+1}_+,\, x_0^adX))$.
    Since $G(X,t) = e^{-|X|^2/(4t)}$, we have for example
    \[
        |\nabla U|^2 G = \left| e^{-|X|^2/(8t)} \nabla U \right|^2,
    \]
    which belongs to $C((0,T);L^1(\R^{d+1}_+,\, x_0^adX))$. The same holds for $U^2G$, $|X\cdot\nabla U|^2G$, and $(\partial_t U)^2G$.

    In particular, Lemma \ref{Lemma convergence hn} also applies to 
    \[
        f(X,t)=\left( X\cdot \nabla U + 2t\partial_t U + \frac{2 s}{p-1}U \right)^2,
    \]
    since
    \[
        g(X,t)=\left( X\cdot \nabla U + 2t\partial_t U + \frac{2 s}{p-1}U \right)e^{-|X|^2/8t}
    \]
    is in $C((0,T);L^2(\R^{d+1}_+,\, x_0^adX))$, and thus $fG = |g|^2 \in C((0,T);L^1(\R^{d+1}_+,\, x_0^adX))$.
\end{Remark}

\begin{Lemma}
    Let $\varepsilon \in (0,1).$ Define $F_n: C([\varepsilon,1]) \mapsto \R$ such that
    $$F_n(f) := n\int_{\varepsilon}^1 t^{n-1} f(t)\, dt,$$
    and
    $$F(f) := f(1).$$
    Then $F_n(f) \to F(f)$.
   \begin{proof}
Clearly $F_n$ is linear for every $n\in \mathbb{N}$. To see that $F_n$ is bounded, we compute
$$\lvert F_n(f) \lvert  \leq n |f|_{\infty} \int_\varepsilon^1 t^ {n-1} \,dt = n\left( \frac{1}{n} -\frac{\varepsilon^n}{n} \right) |f|_{\infty} \leq |f|_{\infty},$$
and thus $\|F_n\| \leq 1$ for every $n \geq 1$.
Since 
$$F_n(t^ k) = n \left(\frac{1}{n+k}- \frac{\varepsilon^{n+k}}{n+k}\right),$$
we see that 
$$\lim_{n \to \infty} F_n(t^ k) =1 = t^ k(1).$$
Since $F_n$ are uniformly bounded and converge to $F$ on the dense subset of polynomials, it follows that  $F_n(f) \to  F(f)$ for every $f \in C([\varepsilon,1]).$
   \end{proof}
\end{Lemma}

\begin{Lemma}\label{lemma epsilon 1}
    Let $\varepsilon \in (0,1)$ and suppose $h_n:[\varepsilon, 1] \to \R$ is a sequence of continuous functions that converges uniformly to  $h:[\varepsilon, 1] \to \R$. Then
    $F_n(h_n) \to F(h).$
    \begin{proof}
      \be  
      \begin{split}
      | F(h) - F_n(h_n) |& \leq | F(h)-F_n(h) | + |F_n(h) - F_n(h_n)|\\
      & \leq |{F(h)-F_n(h)}|\ + \|h-h_n\|_{\infty} \to 0, \\
      \end{split}
      \ee
      since $\|F_n\| \leq 1$ if $n \geq 1$ and $F_n(h) \to F(h)$.
    \end{proof}
\end{Lemma}
\vs

\section{The monotonicity formula}\label{section 5}
We now proceed to derive our main result. First, observe that if $(r,X,\theta) \in \R_+ \times \R^{d+1}_+ \times \mathbb{S}^{n-1}_1,$ $$V_n(z_0,z,y) = V_n \circ F(r,X,\theta) = U(X,r^2/2),$$ and the same applies for the quantities listed in Lemma \ref{chain rule}. Therefore, whenever we are integrating with respect to the $X$ and $r$ variables, the functions will be evaluated in $(X,r^2/2)$ and, if $r=R$ is fixed, the quantities will be evaluated in $(X,R^2/2)=(X,t),$ by \eqref{variables}. Similarly, if we are integrating with respect to the $x$ and $r$ variables, observe that $V(0,z,y) = V\circ F(r,0,x,\theta) = u(x,r^2/2)$. Finally, for any fixed $r=R$ the quantities are evaluated in $(x,R^2/2)=(x,t)$, by \eqref{variables}. We follow this convention by default unless explicit evaluations are provided.

Let 
\be
C_n := \frac{n^{(d+1+a)/2}}{|\mathbb{S}^{n-1}_1|}.
\ee
We have the following proposition.

\begin{Proposition}\label{limit function}
     Let $u \in C^{2s+\varepsilon}_{t,x}(\R^d \times (-T_F,+\infty))$ be a solution of \eqref{back frac} in $[0,T] \subset (-T_F,T_I)$ and suppose its associated extension satisfies $U \in \mathcal{U}([0,T])$. Let $V_n$ be its $n$-dimensional lift.  Define,
\be\label{mathcal E}
\begin{split}
    \mathcal{E}(R) :=& R^{2 s\frac{p+1}{p-1}-d}\frac{1}{2}\int_{\R^{d+1}_+} x_0^a|\nabla U|^2 G (X,R^2/2)\,dX \\
    &\phantom{.}- \frac{\eta_{s}}{p+1}R^{2 s \frac{p+1}{p-1}-d} \int_{\R^d} |u|^{p+1}\tilde{G}(x,R^2/2)\,dx \\
    & \phantom{.}+ R^{2 s \frac{p+1}{p-1}-d-2}\frac{s}{p-1}\int_{\R^{d+1}_+} x_0^a U^2 G(X,R^2/2)\,dX.
\end{split}
\ee
Then, for any $0<\varepsilon<T/2,$
$$C_n\mathcal{E}_n(R) \to \mathcal{E}(R),$$  for every $R \in [\sqrt{2\varepsilon},\sqrt{2(T-\varepsilon)}]$, where $\mathcal{E}_n$ is defined in \eqref{c_ne_n}.
\begin{proof}
    Let $0<\varepsilon<T/2.$ We write
    \be
    C_n \mathcal{E}_n(R) = \mathcal{E}^1_n(R) + \mathcal{E}^2_n(R) + \mathcal{E}^3_n(R),
    \ee
    where,
    \be
    \mathcal{E}^1_n(R) = C_n R^{2 s \frac{p+1}{p-1}-n-d}\frac{1}{2} \int_{\mathbb{R}_{+}^{n+d+1} \cap B^{n+d+1}_R} z_0^a|\nabla V_n|^{2}\, dyd z d z_0,
    \ee
    \be
    \mathcal{E}^2_n(R) = -n^{\frac{1-a}{2}}C_n\frac{\eta_{s}}{p+1}R^{2 s \frac{p+1}{p-1}-n-d} \int_{\partial \mathbb{R}_{+}^{n+d+1} \cap B^{n+d+1}_R}|V_n|^{p+1} \,d zdy,
    \ee
    and,
    \be
    \mathcal{E}^3_n(R) = C_nR^{2 s \frac{p+1}{p-1}-n-d-1} \frac{s}{p-1} \int_{\partial B^{n+d+1}_R  \cap \mathbb{R}_{+}^{n+d+1}} z_0^a V_n^{2}\, d \sigma.
    \ee
    We start by computing the limit of $\mathcal{E}_n^1.$
    By (\ref{preimage ball}) and (\ref{volume in ball}) we have,
    \be
\begin{split}
&C_n R^{-(n-1)}\frac{1}{2}\int_{\mathbb{R}_{+}^{N+1} \cap B^{N+1}_R} z_0^a|\nabla V_n|^{2} \,dyd z d z_0  \\
&= C_n\frac{|\mathbb{S}^{n-1}_1|}{2n^{(d+1)/2}}\int_0^{R} \int_{\R^{d+1}_+} \left(\frac{x_0}{\sqrt{n}}\right)^a\left(n|\nabla U|^2 + 2(X \cdot \nabla U )\partial_t U +r^2 (\partial_tU)^2\right) \left(\frac{r}{R}\right)^{n-1}G_n\,dXdr \\
& = \frac{1}{2}n\int_0^{R} \int_{\R^{d+1}_+} \left(\frac{r}{R}\right)^{n-1}x_0^a|\nabla U|^2 G_n \,dXdr \\
& \phantom{ = .} + \frac{1}{2}\int_0^{R} \int_{\R^{d+1}_+}\left(\frac{r}{R}\right)^{n-1}x_0^a\left(2(X \cdot \nabla U )\partial_t U +r^2 (\partial_tU)^2\right)G_n\,dXdr.
\end{split}
\ee
In particular, we can write,
\be
\mathcal{E}_n^1 (R) = R^{2 s \frac{p+1}{p-1}-d-1} n\int_0^{R} \left(\frac{r}{R}\right)^{n-1}h^{\mbox{\tiny{I}}}_n(r^2/2)\,dr + R^{2 s \frac{p+1}{p-1}-d-1}\int_0^{R} \left(\frac{r}{R}\right)^{n-1}h^{\mbox{\tiny{II}}}_n(r^2/2)\,dr,
\ee
where
$$h^{\mbox{\tiny{I}}}_n(t) :=\frac{1}{2}\int_{\R^{d+1}_+} x_0^a|\nabla U|^2 G_n (X,t)\,dX, $$
and 
$$h^{\mbox{\tiny{II}}}_n(t) :=\int_{\R^{d+1}_+}x_0^a\left(2(X.\nabla U )\partial_t U +2t (\partial_tU)^2\right)G_n (X,t)\,dX.
$$
 Using that $U\in \mathcal{U}([0,T])$ we can apply Lemma \ref{Lemma convergence hn} to show that $h^{\mbox{\tiny{I}}}_n(t)$ converges uniformly to
$$ h^{\mbox{\tiny{I}}}(t) := \frac{1}{2}\int_{\R^{d+1}_+} x_0^a|\nabla U|^2 G (X,t)\,dX,$$
 for $t \in [\varepsilon/2,T-\varepsilon]$. 
 Now, since $R \in [\sqrt{2\varepsilon},\sqrt{2(T-\varepsilon)}]$, 
\be\label{separation terms}
n\int_0^{R} \left(\frac{r}{R}\right)^{n-1}h^{\mbox{\tiny{I}}}_n(r^2/2)\,dr = n\int_{\sqrt{\varepsilon}}^{R} \left(\frac{r}{R}\right)^{n-1}h^{\mbox{\tiny{I}}}_n(r^2/2)\,dr + n\int_0^{\sqrt{\varepsilon}} \left(\frac{r}{R}\right)^{n-1}h^{\mbox{\tiny{I}}}_n(r^2/2)\,dr. 
\ee
 Given that the convergence $h^{\mbox{\tiny{I}}}_n (t)\to h^{\mbox{\tiny{I}}}(t)$ is uniform in $  [\varepsilon/2,T-\varepsilon]$, the convergence $h^{\mbox{\tiny{I}}}_n (R^2/2)\to h^{\mbox{\tiny{I}}}(R^2/2)$ must be uniform for  $R \in [\sqrt{\varepsilon}, \sqrt{2(T-\varepsilon)}]$. We perform the change of variables $\tilde{r}=r/R$, and apply Lemma \ref{lemma epsilon 1} to show,
\be\label{term from e to 1 example}
n\int_{\sqrt{\varepsilon}}^{R} \left(\frac{r}{R}\right)^{n-1}h^{\mbox{\tiny{I}}}_n(r^2/2)\,dr \to R h^{\mbox{\tiny{I}}}(R^2/2) = \frac{R}{2}\int_{\R^{d+1}_+} x_0^a|\nabla U|^2 G (X,R^2/2)\, dX,
\ee
for every $R \in [\sqrt{2\varepsilon},\sqrt{2(T-\varepsilon)}]$. Moreover, since $U\in \mathcal{U}( [0,T])$, we have, 
\be\label{bound h^I}
|h^{\mbox{\tiny{I}}}_n (t) | \leq  \frac{e}{2}\sup_{t \in (0,T)} \int_{\R^{d+1}_+} x_0^a |\nabla U|^2 G (X,t)\,dX < C, 
\ee
for every $t \in
(0,T)$ and some constant $C >0$. 
Therefore,
\be\label{term from 0 to e example}
n\left| \int_0^{\sqrt{\varepsilon}} \left(\frac{r}{R}\right)^{n-1}h^{\mbox{\tiny{I}}}_n(r^2/2)\,dr \right|\leq C R^{1-n} (\sqrt{\varepsilon})^n \leq C\sqrt{2T}\left(\frac{1}{2}\right)^{n/2} \to 0,
\ee
for every $R \in [\sqrt{2\varepsilon},\sqrt{2(T-\varepsilon)}]$, where we used $r^2/2 < T$ for $r \in [0,\sqrt{\varepsilon}]$. 
Combining both computations we obtain,
\be\label{h^i}
n\int_0^{R} \left(\frac{r}{R}\right)^{n-1}h^{\mbox{\tiny{I}}}_n(r^2/2)\,dr \to \frac{R}{2}\int_{\R^{d+1}_+} x_0^a|\nabla U|^2 G(X,R^2/2) \, dX.
\ee
In order to control the term involving $h^{\mbox{\tiny{II}}}_n,$ we use Hölder's inequality and the fact that $U\in \mathcal{U}( [0,T])$, and by a similar argument to the one discussed in (\ref{bound h^I})  we conclude,
\be
|h^{\mbox{\tiny{II}}}_n(t)| < C,
\ee
for some constant $C>0$ and for every $t \in (0,T).$
 Therefore,
\be\label{h^ii}
\begin{split}
\left| \int_0^R \left(\frac{r}{R}\right)^{n-1}h^{\mbox{\tiny{II}}}_n(r^2/2) \,dr \right| &\leq C\int_0^R \left(\frac{r}{R}\right)^{n-1} \,dr \\ & \leq C\frac{R}{n}\\
&\leq C\frac{\sqrt{2T}}{n} \to 0,
\end{split}
\ee
for every $R \in [\sqrt{2\varepsilon},\sqrt{2(T-\varepsilon)}]$.
Combining (\ref{h^i}) and (\ref{h^ii}), we find,
\be\label{E1}
\mathcal{E}^1_n(R) \to  R^{2 s \frac{p+1}{p-1}-d} \frac{1}{2}\int_{\R^{d+1}_+} x_0^a|\nabla U|^2 G (X,R^2/2)\, dX.
\ee
For the second term, $\mathcal{E}^2_n$, we use (\ref{preimage plane}) and (\ref{volume in plane}) to compute,
\be
\begin{split}
-n^{\frac{1-a}{2}}C_n&\frac{\eta_{s}}{p+1}R^{-(n-1)}\int_{\partial \mathbb{R}_{+}^{n+1} \cap B^{N+1}_R}|V_n|^{p+1} \,d zdy  \\ &= -n^{\frac{1-a-d}{2}}C_n\frac{\eta_{s}}{p+1}|\mathbb{S}^{n-1}| \int_0^R \int_{\partial \R^{d+1}_+}|u|^{p+1}\left(\frac{r}{R}\right)^{n-1}\tilde{G}_n\,dxdr \\
& = -n\frac{\eta_{s}}{p+1}\int_0^R \int_{\R^d} |u|^{p+1}\left(\frac{r}{R}\right)^{n-1}\tilde{G}_n\,dxdr.\\ 
\end{split}
\ee
Since $u$ is parabolic Hölder continuous of order $2s+\varepsilon$, by Remark \ref{Lemma for up} we have that
$$\int_{\R^d} |u|^{p+1} \tilde{G}_n(x,t) \,dx \to \int_{\R^d} |u|^{p+1}\tilde{G}(x,t) \,dx$$
uniformly for  $t \in [\varepsilon/2, T-\varepsilon]$
and,
$$\sup_{t \in (0,T)} \int_{\R^d} |u|^{p+1}\tilde{G} \, dx <  +\infty.$$
Then, we proceed as for $h^{\mbox{\tiny{I}}}_n$ in (\ref{separation terms}), (\ref{bound h^I}) and (\ref{term from 0 to e example}) to obtain,
\be\label{E2}
\mathcal{E}_n^2(R) \to - \frac{\eta_{s}}{p+1}R^{2 s \frac{p+1}{p-1}-d} \int_{\R^d} |u|^{p+1}\tilde{G}\,dx,
\ee
for every $R$ in $[\sqrt{2\varepsilon},\sqrt{2(T-\varepsilon)}].$

\vs

For the third term, $\mathcal{E}_n^3,$ we use (\ref{preimage boundary}) and (\ref{volume in H_R}) to deduce,
\be
\begin{split}
\mathcal{E}^3_n(R) &=
C_nR^{2 s \frac{p+1}{p-1}-n-d-1} \frac{s}{p-1} \int_{\partial B^{n+d+1}_R  \cap \mathbb{R}_{+}^{n+d+1}} z_0^a V_n^{2}\, d \sigma,\\
&= R^{2 s \frac{p+1}{p-1}-n-d-1} \frac{s}{p-1}n^{(d+a+1)/2} \int_{\R^{d+1}_+} \left ( \frac{x_0}{\sqrt{n}}\right)^a U^2 n^{-(d+1)/2}R^{n-1}G_n\,dX \\
& = R^{2 s \frac{p+1}{p-1}-d-2} \frac{s}{p-1}\int_{\R^{d+1}_+} x_0^a U^2 G_n\,dX.
\end{split}
\ee
Since $U\in \mathcal{U}([0,T])$, we have  $R^{2 s \frac{p+1}{p-1}-d-2} U^2G (\cdot, R^2/2) 
 \in C([\sqrt{2\varepsilon},\sqrt{2(T-\varepsilon)}]; L^1(R^{d+1}_+, x_0^a dX)), $ and we can directly apply Lemma \ref{Lemma convergence hn} to show,
\be\label{E3}
\mathcal{E}^3_n(R) \to  R^{2 s \frac{p+1}{p-1}-d-2} \frac{s}{p-1} \int_{\R^{d+1}_+} x_0^a U^2G\,dX,
\ee
 for $R \in [\sqrt{2\varepsilon}, \sqrt{2(T-\varepsilon)}].$ 
Combining (\ref{E1}), (\ref{E2}) and (\ref{E3}) the result follows.
\end{proof}
\end{Proposition}

\vs

We now examine the convergence of the derivatives of $C_n\mathcal{E}_n$.

\begin{Proposition}\label{limit derivative}
Let  $U\in \mathcal{U}([0,T])$ and let $0<\varepsilon< T/2$. Define,
\be\label{derivative mathcal E}
\mathcal{D}(R) : = R^{2 s \frac{p+1}{p-1}-d-3}\int_{\R^{d+1}_+} x_0^a\left( X \cdot \nabla U + R^2\partial_t U + \frac{2 s}{p-1}U \right)^2G(X,R^2/2)\,dX.
\ee
Then 
$$\frac{d}{dR}C_n\mathcal{E}_n(R) \to \mathcal{D}(R),$$
 uniformly in $[\sqrt{2\varepsilon},\sqrt{2(T-\varepsilon)}]$. 

 \begin{proof}
Fix  $0<\varepsilon<T/2$.  Using (\ref{der c_n}), we write
     \be
\frac{d}{dR}C_n \mathcal{E}_n(R) = A_n(R) - B_n(R),
\ee
where
\be\label{1st term derivative}
A_n(R) := C_nR^{2 s \frac{p+1}{p-1}-n-d} \int_{\partial B^{n+d+1}_R \cap \mathbb{R}_{+}^{n+d+1}} z_0^{a}\left(\frac{\partial V_n}{\partial r}+\frac{2 s}{p-1} \frac{V_n}{r}\right)^{2} d \sigma,
\ee
and
\be
B_n(R) :=
C_nR^{2s\frac{p+1}{p-1}-1-n-d}
\int_{\mathbb{R}_{+}^{n+d+1} \cap B_{R}^{n+d+1}}
\left(\frac{2s}{p-1}V_n
+(z_0,z,y)\cdot \nabla V_n \right)
(z_0^{a}H \circ \mathcal{F}_n)\,dz_0dzdy.
\ee
To examine the limit of $A_n$, first observe that for any $(z_0,z,y) \in \partial B^{n+d+1}_R \cap \mathbb{R}_{+}^{n+d+1}$, $|(z_0,z,y)|=R$, and therefore
$$\frac{\partial V_n}{\partial r} = \frac{1}{R}(z_0,z,y)\cdot\nabla V_n   = \frac{1}{R}(X\cdot\nabla U + R^2\partial_t U).$$
Then, 
\be
\begin{split}
z_0^{a}\left(\frac{\partial V_n}{\partial r}+\frac{2 s}{p-1} \frac{V_n}{r}\right)^{2} &= \left(\frac{x_0}{\sqrt{n}}\right)^a \left( \frac{1}{R}(X\cdot\nabla U + R^2\partial_t U) + \frac{2 s}{p-1} \frac{U}{R} \right)^2 \\
& = n^{-a/2}\frac{x_0^a}{R^2}\left( X\cdot\nabla U + R^2\partial_t U + \frac{2 s}{p-1}U \right)^2.
\end{split}
\ee
By (\ref{preimage boundary}) and (\ref{volume in H_R}),
\be
\begin{split}
A_n(R)& = C_nR^{2 s \frac{p+1}{p-1}-n-d} \int_{\partial B^{n+d+1}_R \cap \mathbb{R}_{+}^{n+d+1}} z_0^{a}\left(\frac{\partial V_n}{\partial r}+\frac{2 s}{p-1} \frac{V_n}{r}\right)^{2} \,d \sigma,\\
& = R^{2 s \frac{p+1}{p-1}-n-d} \int_{\R^{d+1}_+} \frac{x_0^a}{R^2}\left(X \cdot \nabla U + R^2\partial_t U + \frac{2 s}{p-1}U \right)^2R^{n-1}G_n\,dX.\\
&= R^{2 s \frac{p+1}{p-1}-d-3} \int_{\R^{d+1}_+} x_0^a\left(X\cdot \nabla U + R^2\partial_t U + \frac{2 s}{p-1}U \right)^2G_n\,dX.
\end{split}
\ee
Since $U\in \mathcal{U}([0,T])$,  Remark \ref{lemma for derivative term} and Lemma \ref{Lemma convergence hn}  imply,
$$A_n \to \mathcal{D},$$ 
uniformly for  $R \in [\sqrt{2\varepsilon},\sqrt{2(T-\varepsilon)}].$ 
Finally, we examine the non-homogeneous contribution, and show that it converges uniformly to zero.
First notice that, by (\ref{variables}), we can rewrite
\be
\frac{2s}{p-1}V_n +(z_0,z,y)\cdot \nabla V_n  = \frac{2s}{p-1}U +(x_0,x)\cdot \nabla U +r^2 \partial_tU.
\ee
Then,
\be
\begin{split}
 B_n(R)= R^{2s\frac{p+1}{p-1}-2-d}\int_0^R \int_{\R^{d+1}_+} x_0^aH\left(\frac{2s}{p-1}U+ X\cdot\nabla U + r^2\partial_t U\right)\, \left(\frac{r}{R}\right)^{n-1}G_n\,dXdr.
\end{split}
\ee

Let 
$$s_n(t) := \int_{\R^{d+1}_+} x_0^aH\left(\frac{2s}{p-1}U +X\cdot\nabla U + 2t\partial_t U\right)G_n\,dX.$$
Using that  $U \in \mathcal{U}([0,T])$,
$$|s_n(t)| \leq C,$$
 for $t \in (0,T)$ and some constant $C>0$. Now,
\be
\begin{split}
| B_n(R) |&\leq R^{2s\frac{p+1}{p-1}-2-d}\int_0^R \left(\frac{r}{R}\right)^{n-1}|s_n(r^2/2)|\,dr \\
&\leq R^{2s\frac{p+1}{p-1}-2-d} C\int_0^R\left(\frac{r}{R}\right)^{n-1}\,dr \\& \leq  \frac{C_{\varepsilon,T}}{n},
\end{split}
\ee
for every $R \in [\sqrt{2\varepsilon},\sqrt{2(T-\varepsilon)}]$, where $C_{\varepsilon,T}$ is a constant depending on $T$ and $\varepsilon$. Since the last bound is independent of $R$,  it converges uniformly to zero in $[\sqrt{2\varepsilon},\sqrt{2(T-\varepsilon)}]$.
 \end{proof}
\end{Proposition}

Combining the previous results, we prove our main theorem.

\begin{theorem}\label{main}
   Let  $u \in C^{2s+\varepsilon}_{t,x}(\R^d \times (-T_F,+\infty))$ be a solution of \eqref{back frac}. Let $[0,T] \subset (-T_F,T_I)$ and suppose its associated extension satisfies $U \in \mathcal{U}([0,T])$.  Then, the quantity $($\ref{J}$)$ is non-decreasing in $(0,T)$. Furthermore, its derivative is given by $($\ref{DJ}$)$.

\begin{proof}
Let $0<\varepsilon<T/2,$ and let $t \in [\varepsilon, T-\varepsilon].$ Then, $R \in [\sqrt{2\varepsilon},\sqrt{2(T-\varepsilon)}]$ and, by Proposition \ref{limit derivative} the convergence
$$\frac{d}{dR}C_n\mathcal{E}_n(R) \to \mathcal{D}(R)$$
is uniform in $[\sqrt{2\varepsilon},\sqrt{2(T-\varepsilon)}]$. Furthermore, by Proposition \ref{limit function} we have
\be\label{argument limit}
\lim_{n \to \infty} C_n\mathcal{E}_n(R) = \mathcal{E}(R),
\ee
for every $R \in [\sqrt{2\varepsilon},\sqrt{2(T-\varepsilon)}]$. Then, a standard argument shows,
\be\label{derivative limit}
\frac{d}{dR}\mathcal{E}(R) = \lim_{n \to \infty} \frac{d}{dR} (C_n \mathcal{E}_n)(R) = \mathcal{D}(R),
\ee
for every  $R \in [\sqrt{2\varepsilon},\sqrt{2(T-\varepsilon)}]$. Since  $\varepsilon>0$ is arbitrary, (\ref{argument limit}) and (\ref{derivative limit}) must hold for every $R \in (0,\sqrt{2T}).$

We may now define
$$\mathcal{J}(t) := \frac{1}{(4\pi)^{d/2}\Gamma(s)}\frac{1}{2^{s\frac{p+1}{p-1}-d/2}}\mathcal{E}(\sqrt{2t}),$$
and since $d/dt  = (2t)^{-1/2}d/dR,$
$$\frac{d}{dt}\mathcal{J}(t) = \frac{1}{(4\pi)^{d/2}\Gamma(s)}\frac{1}{2^{s\frac{p+1}{p-1}-d/2}}\frac{1}{\sqrt{2t}}\frac{d\mathcal{E}}{dR} (\sqrt{2t}).$$
Using the expression (\ref{kernel}) for the fundamental solution, and denoting
$$\tilde{\mathcal{G}}(x,t) =\mathcal{G}((0,x),t), $$
we find,
\be
\begin{split}
    \mathcal{J}(t) =& t^{  \frac{2s}{p-1}+1}\left(\frac{1}{2}\int_{\R^{d+1}_+} x_0^{1-2s}|\nabla U|^2 \mathcal{G} \, dX - \frac{\eta_{s}}{p+1} \int_{\R^d} |u|^{p+1}\tilde{\mathcal{G}} \, dx\right) \\
    & + t^{\frac{2s}{p-1}}\frac{s}{2(p-1)}\int_{\R^{d+1}_+} x_0^{1-2s} U^2 \mathcal{G} \, dX,
\end{split}
\ee
and
\be
\begin{split}
\frac{d}{dt}\mathcal{J}(t)
&=
\frac14 t^{ \frac{2s}{p-1}-1}
\int_{\R^{d+1}_+} x_0^{1-2s}
\left( 2t\partial_t U +X \cdot \nabla U
+ \frac{2 s}{p-1}U \right)^2
\mathcal{G}\, dX.
\end{split}
\ee
The result readily follows after rearranging the terms to match those of expressions \eqref{J} and \eqref{DJ} respectively.
\end{proof}
\end{theorem}

\section{Appendix}

Here we prove Lemma \ref{Lemma convergence hn}. We will use the following proposition.
\begin{Proposition}\label{Prop K}
     Let $g \in C([t_0,t_1]; L^1(\R^{d+1}_+, \,x_0^adX)).$
 Then, for every $\varepsilon>0$ there exists a compact set $K \subset \R^{d+1}_+$ such that
 $$\int_{\R^{d+1}_+\setminus K} |g(X,t)| \,x_0^adX < \varepsilon,$$
 for every $t \in [t_0,t_1].$
 \begin{proof} Fix $\varepsilon>0$. Since $g$ is continuous, $g([t_0,t_1]) \subset L^1(\R^{d+1}_+,\, x_0^adX)$ is compact and hence totally bounded. Hence, there exist $g_1,\dots, g_N \in L^1(\R^{d+1}_+,\, x_0^adX)$ such that for any $t \in  [t_0,t_1],$ there is $i \in \{1,\dots,N\}$ with
     $$||g(\cdot,t)-g_i||_{ L^1(\R^{d+1}_+, \,x_0^adX)} < \varepsilon/2.$$
     Since each $g_i$ is integrable,  we can find a compact set $K_i\subset \R^{d+1}_+$ such that
     $$\int_{\R^{d+1}_+ \setminus K_i }|g_i(X)|\,x_0^adX < \varepsilon/2.$$
     Let
     $$K := \bigcup_{i=1}^N K_i,$$
which is compact.  Then, for any $t\in[t_0,t_1]$ and its corresponding $i$,
$$
\int_{\R^{d+1}_+\setminus K} |g(X,t)|\, x_0^a\,dX 
\le \|g(\cdot,t)-g_i\|_{L^1(\R^{d+1}_+,\, x_0^a\,dX)} 
      + \int_{\R^{d+1}_+\setminus K} |g_i(X)|\, x_0^a\,dX 
< \varepsilon.
$$
 \end{proof}
\end{Proposition}

    \begin{proof}[Proof of Lemma \ref{Lemma convergence hn}]
        Let $\varepsilon>0$. First observe that since $G_n \leq e G,$
        $$|h_n(t)| \leq \int_{\R^{d+1}_+}|f|G_n \,x_0^adX \leq e \int_{\R^{d+1}_+} |f|G \,x_0^adX < \infty.$$
        Now,
        \be
        |h_n(t)-h(t)| \leq \int_K |f||G_n-G|\, x_0^adX +\int_{\R^{d+1}_+\setminus K} |f||G_n-G| \, x_0^adX, 
        \ee
        for every compact set $K$.
        Since $fG(\cdot,t) \in  C([t_0,t_1]; L^1(\R^{d+1}_+, x_0^adX))$, by Proposition \ref{Prop K} there exists a compact set $K$ such that
        $$\int_{\R^{d+1}_+ \setminus K} |fG| \, x_0^adX < \frac{\varepsilon}{2(1+e)},$$
        for every $t \in [t_0,t_1].$
        Then, since $G_n \leq eG$
        $$\int_{\R^{d+1}_+\setminus K} |f||G_n-G|\, x_0^adX  \leq (1+e)\int_{\R^{d+1}_+\setminus K} |fG|\, x_0^adX < \frac{\varepsilon}{2}. $$
        Now, 
        $$\int_K |f||G_n-G|\, x_0^adX \leq \max_{(x,t) \in K \times [t_0,t_1]}|f|   \max_{(x,t) \in K \times [t_0,t_1]} |G_n-G| \, \int_K x_0^adX .$$
        Since $K$ is compact and $G_n \to G$ uniformly on $K \times [t_0,t_1]$, there exists $N$ such that, for every  $n\geq N$, the last term is bounded by $\varepsilon/2$.  
Combining the two bounds,  the lemma follows.
        \end{proof}

\bibliographystyle{plain}

\bibliography{bibliography}

\end{document}